\documentclass[11pt]{article}
\usepackage[a4paper,margin=1in]{geometry}
\usepackage[T1]{fontenc}
\usepackage[utf8]{inputenc}
\usepackage{textcomp}
\usepackage{lmodern}
\usepackage{iftex}
\ifPDFTeX
\DeclareUnicodeCharacter{00A7}{\S}
\DeclareUnicodeCharacter{00D7}{\ensuremath{\times}}
\DeclareUnicodeCharacter{2013}{--}
\DeclareUnicodeCharacter{2014}{---}
\DeclareUnicodeCharacter{2018}{`}
\DeclareUnicodeCharacter{2019}{'}
\DeclareUnicodeCharacter{201C}{``}
\DeclareUnicodeCharacter{201D}{''}
\DeclareUnicodeCharacter{2032}{\ensuremath{^\prime}}
\DeclareUnicodeCharacter{2033}{\ensuremath{^{\prime\prime}}}
\DeclareUnicodeCharacter{2192}{\ensuremath{\to}}
\DeclareUnicodeCharacter{2194}{\ensuremath{\leftrightarrow}}
\DeclareUnicodeCharacter{2295}{\ensuremath{\oplus}}
\else
\usepackage{newunicodechar}
\newunicodechar{§}{\S}
\newunicodechar{×}{\ensuremath{\times}}
\newunicodechar{–}{--}
\newunicodechar{—}{---}
\newunicodechar{‘}{`}
\newunicodechar{’}{'}
\newunicodechar{“}{``}
\newunicodechar{”}{''}
\newunicodechar{′}{\ensuremath{^\prime}}
\newunicodechar{″}{\ensuremath{^{\prime\prime}}}
\newunicodechar{→}{\ensuremath{\to}}
\newunicodechar{↔}{\ensuremath{\leftrightarrow}}
\newunicodechar{⊕}{\ensuremath{\oplus}}
\newunicodechar{ğ}{\u{g}}
\newunicodechar{ı}{\i}
\newunicodechar{ş}{\c{s}}
\fi
\usepackage{amsmath,amssymb,amsthm,mathtools,bm}
\usepackage{graphicx}
\usepackage{float}
\usepackage{caption}
\usepackage{booktabs,longtable,array,tabularx}
\usepackage{enumitem}
\usepackage{microtype}
\usepackage{xcolor}
\usepackage[unicode=true]{hyperref}
\usepackage[nameinlink,noabbrev]{cleveref}
\graphicspath{{./}{}}
\providecommand{\tightlist}{%
  \setlength{\itemsep}{0pt}\setlength{\parskip}{0pt}}

\hypersetup{pdftitle={Transpose-odd operator chirality in co-moving multiplicative processes: exact tail-level structure and a detectability obstruction}, pdfauthor={Nihat Çağrı Çalışkan}, colorlinks=true, linkcolor=blue!50!black, citecolor=blue!50!black, urlcolor=blue!50!black}

\theoremstyle{remark}

\title{Transpose-odd operator chirality in co-moving multiplicative processes: exact tail-level structure and a detectability obstruction}
\author{Nihat Çağrı Çalışkan\\[0.45em]\textit{Department of Mathematics, Hacettepe University}\\\textit{06800 Ankara, Türkiye}}
\date{July 15, 2026}
\begin{document}
\maketitle
\begin{abstract}
We ask whether transposition changes the stationary heavy tail of the co-moving recursion \(\mathbf x_{t+1}=Q_tBQ_t^\top\mathbf x_t+Q_t\boldsymbol\eta_t\) in \(d=3\). The model has independent Haar rotations, Gaussian body-frame noise with covariance \(\Sigma\), and an explicit positive metric \(\mathsf g\). The tail index \(\alpha_\star\) depends only on singular values, so \(\alpha_\star(B)=\alpha_\star(B^\top)\).

The structure theorem and polar-twist mechanism are unconditional. In the \(\mathsf g=I\) chart, write \(B=S+\widehat{\boldsymbol\omega}\)\footnotemark:
\[
\Delta_6(B)=-16\det[\boldsymbol\omega,S\boldsymbol\omega,S^2\boldsymbol\omega].
\]
\footnotetext{This degree-six quantity coincides, up to normalization, with the classical sign-indeterminate mixed invariant of a symmetric-plus-skew tensor pair in the \(O(3)\) integrity basis; see Spencer and Rivlin (1962), Smith (1971), and Zheng (1993, 1994). It also coincides, up to normalization, with the turbulence invariant \(\Lambda=\det[\boldsymbol\omega,S\boldsymbol\omega,S^2\boldsymbol\omega]\) of Vigdorovich and Foysi (2015); it agrees with \(-\det[B,B^\top]\) up to a normalization constant. The contribution here is therefore not the algebraic existence of this invariant, but its self-commutator representation \(\Delta_6=-\det[B,B^\top]\), interpreted as ``the determinant of the departure from normality,'' and its identification as intrinsic chirality equivalent to the Kalman-rank/controllability condition for \((S,\boldsymbol\omega)\).}
Then \(B\not\sim_{O(3)}B^\top\), \(\Delta_6(B)\ne0\), \(\operatorname{rank}[\boldsymbol\omega,S\boldsymbol\omega,S^2\boldsymbol\omega]=3\), and controllability of \((S,\boldsymbol\omega)\) are equivalent; intrinsic operator chirality is a Kalman rank condition.

Unlike the exponent, the tail level need not be transpose-invariant. The exact polar-twist identity rotates the anisotropic noise frame under mirroring, showing how level asymmetry can enter. A nonzero asymmetry is established only conditionally and numerically, not by a closed-form example: under Assumption R on the SPD chart, \(\Delta\log C=\langle A,G\rangle+o(\|A\|)\) as \(A\to0\), with \(G=2H\circ[\Sigma,\partial_\Sigma\log C]\); finite-threshold simulations are consistent with this first-order law. Strict-Haar averaging makes every lagged second-order cross-moment transpose-blind; a fourth radial moment retains an odd channel. No finite scalar-weighted radial-moment combination uniformly cancels the even quadratic form while retaining the odd covector.
\end{abstract}

\subsection{1. Introduction}\label{introduction}

Covariance describes dependence but not the direction in which a non-normal operator propagates shocks. We study the part of a body operator \(B\) reversed by the adjoint relative to a positive metric \(\mathsf g\): \(A_{\mathsf g}=(B-B^{\dagger_{\mathsf g}})/2\). In the \(\mathsf g\)-orthonormal chart, this is \(A=(B-B^\top)/2\). We call this transpose-odd content \emph{chirality}: metric-relative operator structure, not spatial-reflection parity.

We use a three-dimensional co-moving Kesten recursion. Independent Haar rotations refresh the direction seen by \(B\); Gaussian innovations have fixed body-frame covariance \(\Sigma\), independent of \(\mathsf g\). The stationary magnitude satisfies \(\mathbb P(\|\mathbf x_\infty\|>z)\sim Cz^{-\alpha_\star}\). The mirror gap is
\[
\Delta\log C:=\log C(B;\Sigma)-\log C(B^\top;\Sigma),
\]
It asks whether transpose-odd propagation moves the tail \emph{level} while the exponent is fixed.

The main results are as follows.

\begin{enumerate}
\item \textbf{Intrinsic chirality is controllability.} In the \(\mathsf g=I\) chart, Theorem 2.3 gives \(\Delta_6(B)=-16\det[\boldsymbol\omega,S\boldsymbol\omega,S^2\boldsymbol\omega]\). Hence \(B\not\sim_{O(3)}B^\top\), \(\Delta_6(B)\ne0\), full Kalman rank, and controllability of \((S,\boldsymbol\omega)\) are equivalent.
\item \textbf{Polar-twist level mechanism.} If \(B=U\Xi V^\top\) and \(O=VU^\top\), then \(B^\top=OBO\) and
\[
C(B^\top;\Sigma)=C(B;O^\top\Sigma O).
\]
Thus mirroring rotates the noise frame exactly. On the positive-definite chart and under Assumption R, this gives
\(\Delta\log C=\langle A,G\rangle+o(\|A\|)\) with \(G=2H\circ[\Sigma,D]\), \(D=\partial_\Sigma\log C\). The confinement theorem separates the commutator, Lyapunov-resolvent, and off-diagonal channels. The indefinite chamber has a separate exact rigidity statement, so the apparent pair-sum pole is not a physical resonance.
\item \textbf{Exponent obstruction.} The fresh-Haar radial reduction makes \(\alpha_\star\) the positive root of \(\mathbb E\|B\mathbf u\|^\alpha=1\), hence a singular-value functional invariant under \(B\mapsto B^\top\). More generally, the moment Lyapunov function of any i.i.d. matrix ensemble is transpose-invariant. Whenever a Kesten index is its root, no transpose-odd signal can enter the exponent.
\item \textbf{Metric-relative low-order channel.} Positive mirrors are parametrized by positive metric rays. Within the stated real-linear unital anti-automorphism class, no \(GL(3)\)-natural metric-free positive transpose exists. At fixed metric, a \(B\)-only odd polynomial scalar first occurs at degree six. An independent covariance opens the complete \(B\)-degree-\(\le2\) family, with common zero criterion \([S_{\mathsf g},\Sigma\mathsf g]=0\).
\item \textbf{Scoped observability limits.} Strict-Haar averaging removes transpose-odd content from every lagged second-order cross-moment, while a fourth radial moment retains an odd term. No finite scalar-weighted radial combination uniformly cancels the even quadratic form over operator directions while retaining the odd covector. The exact metric-ray inverse and common-Gram repair are population/oracle, not sample, identification results.
\end{enumerate}

The fluctuation--dissipation representation of \(D\) requires Assumption UI\(^{\prime\prime}\). Its fixed-lag nullity and lag-window escape are proved separately. The independent-copy coupling proof is given in full in the main text (Lemma 4.8.1). Assumptions R and UI\(^{\prime\prime}\) are not theorems. Supplement S.0 is the controlling logical-status ledger.

\bigskip
\textbf{Related work.} The tail argument uses Kesten--Goldie theory, the asymptotically linear IFS results of Alsmeyer--Brofferio--Buraczewski, stationarity from Bougerol--Picard and Elton--Diaconis--Freedman, and positivity from Buraczewski--Damek. Procesi's orthogonal trace-word theorem supplies the invariant-theory input, and Kubo's fluctuation--dissipation viewpoint frames the response. Transpose-odd propagation differs from asymmetric upper/lower tail dependence (Ang--Chen; Longin--Solnik; Patton) and kinetic-exchange asymmetry.

The nearest work is Troude and Sornette on non-normal random multiplicative processes. Their explicit rotation--reinjection calculation is \(2\times2\), with a separate general-\(N\) approximation; our exact result concerns \(d=3\) with fresh-Haar reinjection. More importantly, the dynamics differ. Their successive operators act on directions produced by earlier steps, preserving eigenvector geometry and alignment. This directional memory lets non-normal products build coherent transient amplification and modify the moment Lyapunov exponent, hence the Kesten tail exponent. Here, fresh independent Haar rotations renew the direction presented to \(B\) and destroy that memory. The radial dynamics becomes a scalar Kesten recursion with i.i.d.\ gains \(\|B\hat{\mathbf u}\|\), switching off multidimensional non-normal amplification. Our model is therefore a strong isotropic-randomization limit of that mechanism, not a reparametrization. Its exponent is transpose-blind because it depends only on the singular values of \(B\); persistent eigenvector geometry can instead modify the exponent in their setting relative to the isotropically randomized benchmark. The two views are complementary: their model retains amplification at the cost of analytic tractability; ours suppresses it to obtain exact structure. Their driving noise is isotropic, so Corollary 4.1.A makes every functional of the norm process mirror-invariant at all orders. The transpose-odd level effect studied here therefore cannot arise in their setting. We use no approximation from their work.

\bigskip
\textbf{Organisation.} Section 2 fixes the model and metric, proves the \(\Delta_6\)/Kalman structure theorem, and gives the calibration. Section 3 proves the exponent result. Section 4 covers the filtered pairing, polar twist, conditional response, coupling, observability, oracle inverse, and open problems. Section 5 states the scope and outlook. Supplement S.A--S.E contains the proof audits and computations.

\subsection{2. Model and the mirror operation}\label{model-and-the-mirror-operation}

\subsubsection{2.1 The co-moving Kesten process}\label{the-co-moving-kesten-process}

In \(d=3\), let
\[
\mathbf x_{t+1}=M_t\mathbf x_t+\boldsymbol\xi_t,
\qquad M_t=Q_tBQ_t^\top,
\qquad \boldsymbol\xi_t=Q_t\boldsymbol\eta_t,
\tag{2.1}
\]
where \(Q_t\stackrel{\rm iid}{\sim}{\rm Haar}(SO(3))\) and \(\boldsymbol\eta_t\stackrel{\rm iid}{\sim}\mathcal N(0,\Sigma)\) are independent. In the body frame \(\mathbf y_t=Q_{t-1}^\top\mathbf x_t\),
\[
\mathbf y_{t+1}=BR_t\mathbf y_t+\boldsymbol\eta_t,
\qquad R_t=Q_t^\top Q_{t-1},
\tag{2.2}
\]
and \(\|\mathbf x_t\|=\|\mathbf y_t\|\) pathwise. Co-rotation fixes the relative orientation of \(B\) and \(\Sigma\). Laboratory-fixed noise defines another model, not a coordinate change.

Fresh Haar makes the direction presented to \(B\) conditionally uniform at each step. Thus the radial gains are i.i.d. with law \(\|B\mathbf u\|\), \(\mathbf u\sim\mathrm{Unif}(S^2)\), and
\[
\mathbb P(\|\mathbf x_\infty\|>z)\sim Cz^{-\alpha_\star},
\qquad
\mathbb E_{\mathbf u}\|B\mathbf u\|^{\alpha_\star}=1.
\tag{2.3}
\]
At the working point, \(\gamma=\mathbb E\log\|B\mathbf u\|=-0.30\), \(\rho(B)=0.770<1<\sigma_{\max}(B)=1.321\), and \(\alpha_\star=6.84\). Section 3 proves existence of the positive limiting constant and the root characterization. Supplements S.A.1--S.A.2 and S.E contain the calibration sweeps, finite-sample exponent diagnostics, and laboratory-fixed-noise control.

\subsubsection{2.2 The symmetric / transpose-odd decomposition}\label{the-symmetric-transpose-odd-decomposition}

In the model chart,
\[
B=S+A,\qquad S=\tfrac12(B+B^\top),\qquad A=\tfrac12(B-B^\top),
\tag{2.4}
\]
so the mirror fixes \(S\) and reverses \(A\): \(B^\top=S-A\). More generally, for a positive metric \(\mathsf g:V\to V^*\),
\[
B^{\dagger_{\mathsf g}}=\mathsf g^{-1}B^\top\mathsf g,\qquad
S_{\mathsf g}=\tfrac12(B+B^{\dagger_{\mathsf g}}),\qquad
A_{\mathsf g}=\tfrac12(B-B^{\dagger_{\mathsf g}}).
\tag{2.5}
\]
Under \(H\in GL(V)\),
\[
B'=HBH^{-1},\qquad \mathsf g'=H^{-\top}\mathsf gH^{-1},\qquad \Sigma'=H\Sigma H^\top,
\]
and \((B')^{\dagger_{\mathsf g'}}=HB^{\dagger_{\mathsf g}}H^{-1}\). The background metric \(\mathsf g:V\to V^*\) and covariance \(\Sigma:V^*\to V\) have opposite tensor variance and are independent inputs; we do not set \(\mathsf g=\Sigma^{-1}\). Supplement S.B.1 gives the full metric bookkeeping.

\bigskip
\textbf{Proposition 2.1 (metric relativity and no canonical positive mirror).}
\emph{Call a real-linear unital involutive algebra anti-automorphism \(\mu\) of \(\operatorname{End}(V)\) positive when \(\operatorname{tr}(\mu(X)X)>0\) for every nonzero \(X\).  In dimension three, positive mirrors are in bijection with positive metric rays,
\[
[\mathsf g]\longleftrightarrow\mu_{\mathsf g},
\qquad
\mu_{\mathsf g}(B)=\mathsf g^{-1}B^\top\mathsf g.
\]
Their moduli space is
\[
GL(3)/\mathbb R^*O(3)\simeq SL(3)/SO(3),
\qquad \dim=5.
\]
Consequently there is no \(GL(3)\)-fixed positive mirror, hence no basis-natural metric-free transpose within this fixed real-linear unital anti-automorphism class.  The well-defined object is instead the equivariant pair-functional \((B,\mathsf g)\mapsto A_{\mathsf g}(B)\).}

\emph{Proof.}
By Skolem--Noether, composing an anti-automorphism with ordinary transpose gives an inner automorphism, so \(\mu(B)=G^{-1}B^\top G\), with \(G\) unique up to nonzero scale.  The involution condition gives \(G^\top=\pm G\).  Odd dimension excludes a nondegenerate skew-symmetric \(G\), while positivity is equivalent to \(\pm G\succ0\).  Thus a positive mirror determines, and is determined by, a positive metric ray.  The \(GL(3)\)-action is transitive and its stabiliser at \(\mathsf g=I\) is \(\mathbb R^*O(3)\); the resulting positive-dimensional homogeneous space has no fixed point. \(\blacksquare\)

\bigskip
\textbf{Proposition 2.2 (fixed-metric scalar threshold).}
\emph{Fix \(\mathsf g\).  Every polynomial scalar \(O(\mathsf g)\)-invariant \(F(B)\) satisfying \(F(B^{\dagger_{\mathsf g}})=-F(B)\) is either zero or has \(B\)-degree at least six.  The bound is sharp in \(d=3\).  One degree-six invariant is
\[
\Delta_6^{\mathsf g}(B):=
\operatorname{tr}\!\left(
B B B^{\dagger_{\mathsf g}} B B^{\dagger_{\mathsf g}} B^{\dagger_{\mathsf g}}
\right)
-
\operatorname{tr}\!\left(
B^{\dagger_{\mathsf g}}B^{\dagger_{\mathsf g}}B B^{\dagger_{\mathsf g}}B B
\right).
\]
In the chart \(\mathsf g=I\), the integer matrix
\[
B_0=\begin{pmatrix}
2&-1&2\\
0&-2&-2\\
-2&-2&-1
\end{pmatrix}
\]
gives \(\Delta_6^I(B_0)=2581-2837=-256\ne0\).}

\emph{Proof.}
Whitening \(\mathsf g\) reduces the claim to simultaneous \(O(3)\)-conjugation of \(B\) and \(B^\top\).  By Procesi's orthogonal invariant theorem (Procesi 1976, Theorem 7.1), polynomial invariants are generated by traces of words in \(B\) and \(B^\top\).  Under the mirror the two letters are interchanged; transposing the resulting word identifies its trace with that of the reversed original word.  Every binary necklace of length at most five is cyclically equivalent to its reversal.  At length six there is exactly one chiral pair, \(\{001011,001101\}\), and the displayed integer witness proves that its trace difference is not identically zero in \(d=3\).  In \(d=2\), writing \(B=S+A\), a reflection in an eigenframe of \(S\) fixes \(S\) and sends \(A\) to \(-A\), so \(B\) and \(B^\top\) are \(O(2)\)-conjugate and every transpose-odd \(B\)-only scalar invariant vanishes.  Thus \(d=3\) is minimal for \(B\)-alone scalar chirality; a pair \((B,\Sigma)\) can already possess a low-degree odd pairing in \(d=2\). \(\blacksquare\)

The degree statement concerns fixed \(\mathsf g\) and \(B\)-degree; it asserts no joint polynomial degree in \((B,\mathsf g)\). It does not exclude nonlinear or data-dependent constructions after relaxing the mirror axioms. It is an \emph{existence} statement, not a claim that degree six detects every chiral operator. Degree six is the first nontrivial odd degree, not a dimension-uniform complete detector; in higher dimensions, supports with odd cycles longer than a triangle may require higher-degree invariants.

The witness has a closed form, turning the lemma into a structure theorem. Write \(A=\widehat{\boldsymbol\omega}\) for the axial vector, \(\widehat{\boldsymbol\omega}\mathbf x=\boldsymbol\omega\times\mathbf x\).

\bigskip
\textbf{Theorem 2.3 (the degree-six invariant is a Krylov determinant).}
\emph{In the chart \(\mathsf g=I\), for every \(B=S+A\) with \(S=S^\top\) and \(A=\widehat{\boldsymbol\omega}\),}
\[
\boxed{\ \Delta_6(B)=-16\,\det\!\big[\boldsymbol\omega,\ S\boldsymbol\omega,\ S^2\boldsymbol\omega\big].\ }
\tag{2.3a}
\]

\emph{Proof.} Both sides are polynomials of bidegree \((3,3)\) in \((S,A)\); the finite verification is given in Supplement S.B.1.  The two necklaces share the prefix \(BBB^\top\) and differ by one transposition, so
\[
\Delta_6=\operatorname{tr}\!\big(BBB^\top[B,B^\top]B^\top\big),
\qquad [B,B^\top]=[S+A,S-A]=2[A,S].
\]
Expanding the four remaining slots leaves twelve bidegree-\((3,3)\) words: eight equal \(-t\) and four vanish, where \(t:=\operatorname{tr}(S^2ASA^2)\), hence \(\Delta_6=-16t\).  On the other hand, \(\boldsymbol\omega\times S\boldsymbol\omega=AS\boldsymbol\omega\), \(\boldsymbol\omega\boldsymbol\omega^\top=A^2+\|\boldsymbol\omega\|^2I\), and \(\operatorname{tr}(S^3A)=0\) give
\[
\det[\boldsymbol\omega,S\boldsymbol\omega,S^2\boldsymbol\omega]
=-\operatorname{tr}(SAS^2A^2)=t.
\]
Combining the reductions proves (2.3a). \(\blacksquare\)

The constant is derived, not fitted: the \(2\) is the commutator normalisation and the \(-8\) the signed word count. At \(B_0\), (2.3a) returns \(-256=-16\cdot16\), so the integer witness reads the identity at one point. Symbolically, \(S=S^\top\) has six free entries and \(\boldsymbol\omega\) three. Expanding both sides over these nine indeterminates gives \(\Delta_6(B)+16\det[\boldsymbol\omega,S\boldsymbol\omega,S^2\boldsymbol\omega]\equiv0\), equivalently \(\Delta_6=-16\operatorname{tr}(S^2ASA^2)\). Supplement~S.B.1 gives the full finite word expansion without numerical evaluation.

In the eigenframe of a simple-spectrum \(S\), with \(\boldsymbol\omega'\) the axial vector expressed there,
\[
\det[\boldsymbol\omega,S\boldsymbol\omega,S^2\boldsymbol\omega]
=\omega'_1\omega'_2\omega'_3\prod_{i<j}(s_j-s_i),
\tag{2.3b}
\]
a Vandermonde in the spectrum times the axial-component product. The vector \(\boldsymbol\omega\) is \emph{cyclic} for \(S\), equivalently the single-input pair \((S,\boldsymbol\omega)\) --- the control system \(\dot{\mathbf x}=S\mathbf x+\boldsymbol\omega u\) --- is \emph{controllable}, exactly when its Kalman matrix \([\boldsymbol\omega,S\boldsymbol\omega,S^2\boldsymbol\omega]\) is nonsingular.

\bigskip
\textbf{Corollary 2.4 (chirality in isolation is controllability).}
\emph{In the \(\mathsf g=I\) chart --- equivalently, after \(\mathsf g\)-whitening, with orthogonal similarity understood as \(O(\mathsf g)\)-similarity before whitening --- the following are equivalent for every \(B=S+A\) in \(d=3\), with no hypothesis on the spectrum of \(S\):}
\begin{enumerate}
\item[(i)] \(B\) \emph{is orthogonally similar to} \(B^\top\);
\item[(ii)] \(\Delta_6(B)=0\);
\item[(iii)] \(\det[\boldsymbol\omega,S\boldsymbol\omega,S^2\boldsymbol\omega]=0\);
\item[(iv)] \(\boldsymbol\omega\) \emph{is not a cyclic vector for} \(S\), \emph{equivalently the pair} \((S,\boldsymbol\omega)\) \emph{is not controllable}.
\end{enumerate}

\emph{Proof.} (ii)\(\Leftrightarrow\)(iii) is (2.3a) and (iii)\(\Leftrightarrow\)(iv) is the Kalman rank condition.  For (i)\(\Rightarrow\)(ii), \(\Delta_6\) is an \(O(3)\)-conjugation invariant and \(\Delta_6(B^\top)=-\Delta_6(B)\), so \(B^\top=PBP^\top\) forces \(\Delta_6(B)=0\).  For (iii)\(\Rightarrow\)(i), suppose first that \(S\) has simple spectrum.  Any \(P\in O(3)\) satisfying \(PSP^\top=S\) is, in the eigenframe, a sign involution \(\operatorname{diag}(\varepsilon)\), while \(PAP^\top=-A\) requires \(\varepsilon_i\varepsilon_j=-1\) whenever \(A'_{ij}\ne0\).  By (2.3b), (iii) forces some \(\omega'_k=0\), which removes one of the three off-diagonal entries of \(A'\); the two remaining constraints are satisfiable and give \(PBP^\top=B^\top\).  If \(S\) has a repeated eigenvalue, its minimal polynomial has degree \(\le2\), hence the Krylov matrix is singular; Lemma 2.4.A supplies \(P\) directly. \(\blacksquare\)

\bigskip
\textbf{Lemma 2.4.A (repeated spectrum).}
\emph{If \(S\) has an eigenvalue of multiplicity at least two, then \(B\) is orthogonally similar to \(B^\top\) for every \(A\).}

\emph{Proof.} Let \(E\) be an eigenspace with \(\dim E\ge2\) and decompose \(\boldsymbol\omega=\boldsymbol\omega_E+\boldsymbol\omega_{E^\perp}\).  Take \(P\) to be the identity on every other eigenspace and, on \(E\), the reflection through a hyperplane of \(E\) containing \(\boldsymbol\omega_E\).  Then \(PSP^\top=S\), \(P\boldsymbol\omega=\boldsymbol\omega\), and \(\det P=-1\), so \(PAP^\top=\det(P)\widehat{P\boldsymbol\omega}=-A\) and \(PBP^\top=B^\top\). \(\blacksquare\)

Corollary 2.4 gives the algebraic form of our thesis. The core is chiral in isolation exactly when \((S,\boldsymbol\omega)\) is controllable. Otherwise every \(O(3)\)-conjugation-invariant \(B\)-only transpose-odd scalar vanishes, and misalignment of \(\Sigma\) must carry the mirror effect. The canonical operator of \S2.4 satisfies (iii), isolating the covariance-orientation mechanism rather than a construction defect.

\subsubsection{2.3 The mirror observable}\label{the-mirror-observable}

Corollary 3.4 gives \(\alpha_\star(B)=\alpha_\star(B^\top)\). The transpose-odd observable is therefore the finite-depth level contrast
\[
dq_p=\log q_p(\|\mathbf x^{(B)}\|)-\log q_p(\|\mathbf x^{(B^\top)}\|).
\tag{2.6}
\]
For two tails with common exponent \(\alpha_\star\), \(dq_p\to\alpha_\star^{-1}\log[C(B;\Sigma)/C(B^\top;\Sigma)]\) along the tail. This ordered-pair contrast changes sign when exchanging \(B\) and \(B^\top\).

\subsubsection{2.4 Calibration and working point}\label{calibration-and-working-point}

Use the canonical representative
\[
S=\operatorname{diag}(-c,0,c),\qquad c=0.849181,\qquad
A=A(\boldsymbol\omega),\quad\boldsymbol\omega=(0.451,0.387,0),
\tag{2.7}
\]
with scale fixed by \(\gamma=-0.30\). It has \(\sigma(B)=(1.321,0.597,0.219)\), \(\rho(B)=0.770\), and \(\alpha_\star\simeq6.84\). Equal spacing is a calibration convention, not a theorem hypothesis. Setting \(\omega_3=0\) gives
\[
PBP=B^\top,\qquad P=\operatorname{diag}(1,1,-1).
\tag{2.8}
\]
Corollary 2.4 applies directly: \(S\) is diagonal, so \(\boldsymbol\omega'=\boldsymbol\omega\), and \(\omega_3=0\) forces \(\det[\boldsymbol\omega,S\boldsymbol\omega,S^2\boldsymbol\omega]=0\) by (2.3b). The pair \((S,\boldsymbol\omega)\) is uncontrollable, with (2.8) as its sign involution. Thus the canonical core is \emph{achiral in isolation}: every \(B\)-only transpose-odd scalar vanishes, and the failure of \(\Sigma\) to share the symmetry carries the entire mirror effect.

Consequently \(\operatorname{law}\|\mathbf x^{(B^\top)}\|_\Sigma=\operatorname{law}\|\mathbf x^{(B)}\|_{P\Sigma P}\). The mirror gap therefore vanishes exactly, at every order and without a regularity hypothesis, on the \emph{symmetry-enforced null locus}
\[
P\Sigma P=\Sigma\iff\Sigma_{13}=\Sigma_{23}=0.
\tag{2.9}
\]
This codimension-two locus is strictly wider than isotropy and the commuting locus \([S,\Sigma]=0\). Here the latter means diagonal \(\Sigma\) and has codimension three: \(\Sigma=\bigl(\begin{smallmatrix}2&0.6&0\\0.6&1.5&0\\0&0&1.2\end{smallmatrix}\bigr)\) does not commute with \(S\) but satisfies (2.9). Thus Theorem 4.9's commuting hypothesis is sufficient but not sharp here. We call (2.9) symmetry-enforced, not a characterisation of the full zero set; accidental zeros elsewhere remain possible. The locus \(I_1=0\) kills the degree-\((2,1)\) filtered channel and one-step fourth-moment odd term. It is not generally the null locus of the first derivative of \(\Delta\log C\). Supplements S.A.2 and S.E.2 contain the deterministic quadrature, calibration plot, and checkpoint audit.

\subsection{3. Exponent invariance: radial reduction and tail index}\label{exponent-invariance-the-radial-reduction-and-the-tail-index}

\subsubsection{3.1 Radial reduction}\label{the-radial-process-is-markov}

\bigskip
\textbf{Lemma 3.1 (radial reduction).} \emph{The radius is pathwise an iterated function system of i.i.d. random Lipschitz maps:}
\[
r_{n+1}=\Psi_n(r_n):=\|r_nB\mathbf V_n+\boldsymbol\eta_n\|,
\tag{3.1}
\]
\emph{where \(\mathbf V_n=Q_n^\top\widehat{\mathbf x}_n\stackrel{\rm iid}{\sim}\mathrm{Unif}(S^2)\), \(\boldsymbol\eta_n\stackrel{\rm iid}{\sim}\mathcal N(0,\Sigma)\), and each pair is independent of the past.}

\emph{Proof.} Orthogonality gives \(r_{n+1}=\|r_nB Q_n^\top\widehat{\mathbf x}_n+\boldsymbol\eta_n\|\). Conditional on the past, Haar invariance makes \(Q_n^\top\widehat{\mathbf x}_n\) uniform and independent of \(\boldsymbol\eta_n\), with a law not depending on the past. Iterating this conditional product law makes the driving pairs i.i.d. and proves the Markov/IFS statement. \(\square\)

For the noiseless product \(\Pi_n=M_n\cdots M_1\), each column radius is a product of i.i.d. gains \(A_j=\|B\mathbf V_j\|\). Comparing the operator norm with the three columns gives
\[
\gamma_{\rm top}=\mathbb E\log A=\gamma,
\qquad
\Lambda(s):=\lim_{n\to\infty}\frac1n\log\mathbb E\|\Pi_n\|^s
=\log M(s),\qquad M(s):=\mathbb EA^s.
\tag{3.2}
\]

\subsubsection{3.2 Tail theorem}\label{the-tail-index}

The radial maps are asymptotically linear:
\[
|\Psi_n(r)-A_nr|\le\|\boldsymbol\eta_n\|,
\qquad A_n=\|B\mathbf V_n\|.
\tag{3.3}
\]

\bigskip
\textbf{Lemma 3.2 (a-priori moments).} \emph{For every \(s\in(0,\alpha_\star)\), \(\mathbb E r_\infty^s<\infty\).}

\emph{Proof.} Log-convexity and \(M(0)=M(\alpha_\star)=1\) give \(M(s)<1\) inside the interval. From (3.3), Minkowski for \(s\ge1\), and subadditivity for \(s<1\), the finite-time iterates started at zero have uniformly bounded \(s\)-moments. Mean contraction \(\mathbb E\log A<0\) gives convergence to the unique stationary radius, and Fatou closes the bound. \(\square\)

\bigskip
\textbf{Theorem 3.3 (tail index and positive level).} \emph{Let \(B\in GL(3)\), \(\Sigma\succ0\), and \(A=\|B\mathbf V\|\) for \(\mathbf V\sim\mathrm{Unif}(S^2)\). Assume \(\mathbb E\log A<0\) and \(\mathbb EA^{\alpha_\star}=1\) for some \(\alpha_\star>0\). Then the vector recursion has a unique stationary law and}
\[
\lim_{z\to\infty}z^{\alpha_\star}\mathbb P(\|\mathbf x_\infty\|>z)=C\in(0,\infty),
\qquad M(\alpha_\star)=1.
\]

\emph{Proof.} \emph{Existence.} By (3.1) the radial chain is an iterated function system of random Lipschitz maps: the reverse triangle inequality gives \(|\Psi(r)-\Psi(r')|\le\|B\mathbf V\|\,|r-r'|\), so \(L(\Psi)=A\), with equality as \(r\to\infty\). Hence \(\mathbb E\log L(\Psi)=\gamma<0\), while \(\mathbb E\log^+|\Psi(0)|=\mathbb E\log^+\|\boldsymbol\eta\|<\infty\) because \(\boldsymbol\eta\) is Gaussian. Goldie's Theorem 2.1 gives a unique stationary radial law with \(R\stackrel d=\Psi(R)\). The backward iteration converges almost surely to \(R\), which is measurable with respect to past innovations and independent of the current \((\mathbf V,\boldsymbol\eta)\); this independence is used in Lemma 3.2 and below. The vector law follows from \(\gamma_{\rm top}=\gamma<0\) by (3.2), through the Bougerol--Picard criterion.

\emph{Goldie's hypotheses.} The moment function \(M(s)=\mathbb EA^s\) is finite for every \(s\ge0\) because \(A\le\sigma_{\max}\), and \(M(0)=M(\alpha_\star)=1\); Goldie's Lemma 2.2 yields \(0<m:=\mathbb E[A^{\alpha_\star}\log A]<\infty\). The law of \(A\) has a density on \([\sigma_{\min},\sigma_{\max}]\), so \(\log A\) is non-arithmetic. For the perturbation condition (2.16), first suppose \(0<\alpha_\star\le1\). Subadditivity of \(x\mapsto x^{\alpha_\star}\) and (3.3) give
\[
\mathbb E\big|\Psi(R)^{\alpha_\star}-(AR)^{\alpha_\star}\big|
\le \mathbb E\|\boldsymbol\eta\|^{\alpha_\star}<\infty.
\]
If \(\alpha_\star>1\), the mean-value inequality, Lemma 3.2, and Gaussian moments instead give
\[
\mathbb E\big|\Psi(R)^{\alpha_\star}-(AR)^{\alpha_\star}\big|
\le c_{\alpha_\star}\,\mathbb E\!\left[(AR)^{\alpha_\star-1}\|\boldsymbol\eta\|+\|\boldsymbol\eta\|^{\alpha_\star}\right]<\infty.
\]
\emph{Finiteness.} Goldie's Corollary 2.4 applies to the measurable random map \(\Psi\), the fixed-point law \(R\stackrel d=\Psi(R)\), the independence of \(R\) from \((\Psi,A)\), and the perturbation estimate above. It gives \(\lim_{z\to\infty}z^{\alpha_\star}\mathbb P(R>z)=C\), with \(C=N(\alpha_\star)/(\alpha_\star m)\) by his (2.18), where \(N(\alpha_\star)=\mathbb E[\Psi(R)^{\alpha_\star}-(AR)^{\alpha_\star}]\). This step gives only \(C\ge0\).

\emph{Positivity.} The reverse triangle inequality gives the pathwise minorant \(\Psi(r)\ge Ar+Q\), with the same multiplier \(A=\|B\mathbf V\|\ge0\) and \(Q=-\|\boldsymbol\eta\|\). This is hypothesis (2.2) of Buraczewski and Damek (2017). Their Theorem 1.3 applies because \(\mathbb E\log A<0\), \(\mathbb EA^{\alpha_\star}=1\), \(\mathbb E[|Q|^{\alpha_\star}+A^{\alpha_\star}\log^+A]<\infty\), \(\log A\) is non-arithmetic, and \(\mathbb P(Ax+Q=x)<1\) for every \(x\), the last fact following from the continuous law of \(Q\). Their Lemma 2.6 additionally requires \(\mathbb E\log L(\Psi)<0\), \(\mathbb E\log^+|\Psi(0)|<\infty\), and an unbounded-above stationary law. The first two were established above; the third follows from Anderson's inequality, since \(\mathbb P(\Psi(r)>t\mid r,\mathbf V)\ge\mathbb P(\|\boldsymbol\eta\|>t)>0\). Lemma 2.6 therefore gives \(\varepsilon>0\) with \(\mathbb P(R>t)\ge\varepsilon t^{-\alpha_\star}\) for large \(t\), and hence \(C\ge\varepsilon>0\). No restriction \(\alpha_\star\ge1\) or \(\alpha_\star\ge2\) enters. The full audit is Supplement S.A.3--S.A.4. \(\square\)

\subsubsection{3.3 Exponent invariance}\label{exponent-invariance-under-the-mirror}

\bigskip
\textbf{Corollary 3.4 (singular-value dependence).} \emph{\(\alpha_\star\) depends on \(B\) only through its singular values; thus \(\alpha_\star(B)=\alpha_\star(B^\top)=\alpha_\star(UBW)\) for orthogonal \(U,W\).}

\emph{Proof.} The law of \(\|B\mathbf u\|^2=\mathbf u^\top B^\top B\mathbf u\) for uniform \(\mathbf u\) depends only on the eigenvalues of \(B^\top B\). These are the squared singular values and agree with those of \(BB^\top\), the Gram matrix for \(B^\top\). The root in Theorem 3.3 is therefore unchanged. \(\square\)

\bigskip
\textbf{Proposition 3.5 (generic i.i.d. transpose boundary).} \emph{Let \(M_n\) be i.i.d. matrices and let \(\Lambda(s)=\lim_n n^{-1}\log\mathbb E\|M_n\cdots M_1\|^s\) exist and be finite. For every transpose-invariant submultiplicative norm, \(\Lambda_M(s)=\Lambda_{M^\top}(s)\). Consequently, whenever an affine-recursion tail index is characterized by \(\Lambda(\alpha)=0\), it is invariant under transposition.}

\emph{Proof.} Since \((M_n\cdots M_1)^\top=M_1^\top\cdots M_n^\top\), transpose invariance of the norm and exchangeability of the i.i.d. tuple give the exact finite-\(n\) equality
\(\mathbb E\|M_n\cdots M_1\|^s=\mathbb E\|M_n^\top\cdots M_1^\top\|^s\). Taking the limit proves the claim. \(\square\)

Fresh-Haar isotropy strengthens the generic statement to singular-value-only dependence. Without isotropy, the projective eigenmeasure may retain other transpose-even information. Temporal dependence leaves the i.i.d. boundary and can move exponents, but lies outside this model. Supplements S.A.5--S.A.6 and S.E.2 give the full boundary taxonomy and numerical checks.

\subsection{4. Transpose-odd structure of the tail level: polar-twist mechanism, conditional first-order response, and lab-frame identifiability}\label{the-chiral-level-effect-polar-twist-mechanism-first-order-response-and-lab-frame-identifiability}

\emph{Notation.} Unless noted, \(\mathbf x_t\) is the lab state and \(\mathbf y_t\) the body state; \(B=S+A\), with \(S=S^\top\), \(A=-A^\top\), and \(\boldsymbol\eta_t\sim N(0,\Sigma)\). We use \(\langle X,Y\rangle=\operatorname{tr}(X^\top Y)\) in the \(\mathsf g=I\) chart. The metric-relative objects are \(B^{\dagger_{\mathsf g}},S_{\mathsf g},A_{\mathsf g}\) and \(\mathsf R=\Sigma\mathsf g\). Section 3's moment function is \(M\); \(\Xi\) is the SVD singular-value diagonal. Supplement S.E.2 indexes checkpoints C1--C31.

\subsubsection{4.0 What §4 establishes}\label{what-4-establishes}

Section 4 identifies the exponent-external observable and its mechanism. The commutator \([S_{\mathsf g},\Sigma\mathsf g]\) filters the metric-relative odd channel (Theorem 4.0). The exact polar twist turns transposition into a noise-frame rotation (Theorem 4.1). Under Assumption R, the first-order level response is its Lyapunov-resolvent image (Theorems 4.5 and 4.7); §4.2.4 gives the fluctuation--dissipation representation. The strict-Haar observability results are logically separate. The exact inverse result of §4.4 is population/oracle only.

\bigskip
\textbf{Theorem 4.0 (metric-relative filtered pairing).}
\emph{Let \(\mathsf g:V\to V^*\) be a positive metric, \(\Sigma:V^*\to V\) a positive covariance/reference tensor, and \(B:V\to V\), with \(\dim V=3\). Define
\[
\mathsf R:=\Sigma\mathsf g,
\qquad
S_{\mathsf g}:=\tfrac12(B+B^{\dagger_{\mathsf g}}),
\qquad
A_{\mathsf g}:=\tfrac12(B-B^{\dagger_{\mathsf g}}).
\tag{4.0a}
\]
Then \(\mathsf R'=H\mathsf RH^{-1}\) under the basis change of §2.2, and \(\mathsf R^{\dagger_{\mathsf g}}=\mathsf R\). For integers \(a>b\ge0\), put
\[
\mathcal D^{\mathsf g}_{a,b}:=
\operatorname{tr}\!\left(B\mathsf R^aB^{\dagger_{\mathsf g}}\mathsf R^b\right)
-\operatorname{tr}\!\left(B^{\dagger_{\mathsf g}}\mathsf R^aB\mathsf R^b\right).
\]
Each \(\mathcal D^{\mathsf g}_{a,b}\) is a \(GL(V)\)-invariant scalar of the triple \((B,\mathsf g,\Sigma)\), is odd under the single fixed mirror \(B\mapsto B^{\dagger_{\mathsf g}}\), and factorises as
\[
\boxed{\mathcal D^{\mathsf g}_{a,b}
=2\operatorname{tr}\!\left[
A_{\mathsf g}\bigl(\mathsf R^aS_{\mathsf g}\mathsf R^b
-\mathsf R^bS_{\mathsf g}\mathsf R^a\bigr)\right].}
\tag{4.0b}
\]
Let \(\mathcal P^-_{\mathsf g,\le2}\) denote the mirror-odd part, of \(B\)-degree at most two, of the polynomial \(O(\mathsf g)\)-invariants of \((B,\mathsf R)\), and let
\[
\mathcal E_{\mathsf R}:=\mathbb R[e_1(\mathsf R),e_2(\mathsf R),e_3(\mathsf R)]
\]
be the \(B\)-degree-zero coefficient ring. Then
\[
\mathcal P^-_{\mathsf g,\le2}
=\mathcal E_{\mathsf R}\mathcal D^{\mathsf g}_{1,0}
\oplus\mathcal E_{\mathsf R}\mathcal D^{\mathsf g}_{2,0}
\oplus\mathcal E_{\mathsf R}\mathcal D^{\mathsf g}_{2,1}.
\tag{4.0c}
\]
Finally, for fixed \(S_{\mathsf g}\) and \(\mathsf R\), all these filtered odd invariants vanish for every \(\mathsf g\)-skew \(A_{\mathsf g}\) if and only if
\[
[S_{\mathsf g},\mathsf R]=0.
\tag{4.0d}
\]
In the paper's \(\mathsf g=I\) chart, \(\mathsf R=\Sigma\) and \(\mathcal D^I_{1,0}=2\langle A,[S,\Sigma]\rangle=2I_1\). Thus the independent covariance lowers the first available odd scalar channel from \(B\)-degree six to bidegree \((2,1)\).}

\emph{Proof.}
The transformation laws in §2.2 and trace cyclicity give scalar invariance and mirror parity; expanding \(B=S_{\mathsf g}+A_{\mathsf g}\) gives (4.0b). Whiten \(\mathsf g=K^\top K\). Then \(KBK^{-1}\) is acted on by ordinary transpose and \(K\Sigma K^\top\) is symmetric positive-definite. Procesi's orthogonal trace-word theorem (Procesi 1976, Theorem 7.1) generates the polynomial invariants by traces of words in these matrices and the transpose of the first; restriction to the \(\mathsf g\)-self-adjoint \(\mathsf R\)-subspace is surjective on invariants by compact \(O(\mathsf g)\) Reynolds averaging. Words of \(B\)-degree zero or one are mirror-even. At \(B\)-degree two, the only nonzero mirror-odd necklace differences arise from one \(B\) and one \(B^\top\), separated by two \(\mathsf R\)-arcs. Cayley--Hamilton reduces the arc powers to \(0,1,2\) over \(\mathcal E_{\mathsf R}\), leaving precisely \((a,b)=(1,0),(2,0),(2,1)\). To see that the sum is direct, take \(\mathsf R=\operatorname{diag}(\lambda_1,\lambda_2,\lambda_3)\) with simple spectrum. On the three edge variables \(A_{ij}S_{ij}\), the coefficient rows of the displayed generators are
\[
4(\lambda_j-\lambda_i)
\bigl(1,\lambda_i+\lambda_j,\lambda_i\lambda_j\bigr),
\qquad 1\le i<j\le3.
\]
Their determinant is, up to sign, \(4^3\prod_{i<j}(\lambda_i-\lambda_j)^2\), so no nontrivial \(\mathcal E_{\mathsf R}\)-linear relation exists; genericity and polynomial density complete the argument. If \([S_{\mathsf g},\mathsf R]=0\), (4.0b) vanishes. Conversely \(\mathcal D^{\mathsf g}_{1,0}=2\operatorname{tr}(A_{\mathsf g}[\mathsf R,S_{\mathsf g}])\), and the trace pairing on the \(\mathsf g\)-skew subspace is nondegenerate, so vanishing for every \(A_{\mathsf g}\) forces the commutator to vanish. \(\blacksquare\)

\emph{Metric-ray and whitening caution.}
Replacing \(\mathsf g\) by \(c\mathsf g\), \(c>0\), leaves the mirror unchanged but sends \(\mathsf R\mapsto c\mathsf R\) and \(\mathcal D^{\mathsf g}_{a,b}\mapsto c^{a+b}\mathcal D^{\mathsf g}_{a,b}\). Thus zeros, parity, and (4.0d) are ray-independent. A normalized magnitude can use
\[
\overline{\mathsf R}:=\frac{\mathsf R}{\det(\mathsf R)^{1/3}},
\qquad
\overline{\mathcal D}^{\mathsf g}_{a,b}:=
\frac{\mathcal D^{\mathsf g}_{a,b}}{\det(\mathsf R)^{(a+b)/3}}.
\]
Imposing \(\mathsf g=\Sigma^{-1}\) is not passive whitening: it changes the mirror to \(B\mapsto\Sigma B^\top\Sigma^{-1}\), makes \(\mathsf R=I\), and annihilates all three low-order channels. Starting from \(\mathsf g=I\), passive covariance whitening with \(H=\Sigma^{-1/2}\) gives \((\Sigma',\mathsf g')=(I,\Sigma)\), not \((I,I)\), and preserves \(\mathsf R'\sim\mathsf R\). Resetting the whitened metric to \(I\) is a new model, equivalent in the original coordinates to \(\mathsf g=\Sigma^{-1}\). We instead keep the independent background metric \(\mathsf g=I\) and body-frame innovation covariance \(\Sigma\).

\begin{center}\rule{0.5\linewidth}{0.5pt}\end{center}

\subsubsection{4.1 The polar-twist mechanism and its chamber structure}\label{the-polar-twist-mechanism-and-its-chamber-structure}

\paragraph{4.1.1 Transpose as a noise-frame rotation}\label{transpose-as-a-noise-frame-rotation}

Write the singular value decomposition \(B=U\Xi V^\top\), where \(\Xi\) is the singular-value diagonal and \(D\) remains reserved for the noise gradient \(\partial_\Sigma\log C\). Define the \textbf{twist} \(O:=VU^\top\in O(3)\), the transpose of the orthogonal polar factor of \(B\). The following identity is elementary and exact.

\bigskip
\textbf{Theorem 4.1 (Polar-twist).} \emph{For every invertible \(B\) and every unit vector \(\mathbf u\), \(B^\top\mathbf u = O\,B\,(O\mathbf u)\), where \(O=VU^\top\) is the transpose of the orthogonal polar factor of \(B\) (unique for invertible \(B\)). Consequently the tail constant of the fresh-Haar norm chain satisfies}
\[
\begin{aligned}
C(B^\top;\Sigma)&=C\!\left(B;\,O^\top\Sigma O\right),\\
\Delta\log C
&:=\log C(B;\Sigma)-\log C(B^\top;\Sigma)\\
&=\log C(B;\Sigma)-\log C(B;O^\top\Sigma O).
\end{aligned}
\]

\emph{Proof.} From \(B=U\Xi V^\top\), \(B^\top=V\Xi U^\top=(VU^\top)(U\Xi U^\top)=O\,(U\Xi U^\top)\); and \(U\Xi U^\top=U\Xi V^\top\,VU^\top=B\,O\). Hence \(B^\top=O B O\), i.e.~\(B^\top\mathbf u=OB(O\mathbf u)\). The norm chain \(r_{t+1}=\|r_tB^\top\hat{\mathbf u}_t+\boldsymbol\eta_t\|\) driven by \((\hat{\mathbf u}_t,\boldsymbol\eta_t)\) is therefore pathwise identical, after the deterministic relabelling \(\hat{\mathbf u}\mapsto O\hat{\mathbf u}\), \(\boldsymbol\eta\mapsto O^\top\boldsymbol\eta\), to the chain with operator \(B\) and noise \(O^\top\boldsymbol\eta\sim\mathcal N(0,O^\top\Sigma O)\); the fresh-Haar direction law is \(O\)-invariant. The tail constants therefore coincide. \(\blacksquare\)

\emph{Remark.} The standing Goldie / Buraczewski--Damek conditions of §3 give existence, uniqueness, and positivity of \(C\) for both \((B,\Sigma)\) and \((B,O^\top\Sigma O)\); we invoke rather than re-prove them. The elementary exact identity differs from standard operator-transpose duality in multivariate Kesten theory. Here, fresh Haar makes mirroring the \emph{noise-frame} rotation \(\Sigma\mapsto O^\top\Sigma O\).

Transposition acts by \emph{rotating the noise frame} through the relational angle \(O\). Chirality is invisible against isotropic noise and can appear only through the misalignment \([O,\Sigma]\).

Four corollaries give the consequences.

\bigskip
\textbf{Corollary 4.1.A (isotropy hard-null).} \emph{If \([O,\Sigma]=0\) --- in particular if \(\Sigma\propto I\) in the \(\mathsf g=I\) chart --- then \(O^\top\Sigma O=\Sigma\) and the entire law of the norm process is mirror-invariant, so every functional of \(\{\|\mathbf x_t\|\}\) (level, quantile, spectrum) carries exactly zero chiral signal, at all orders.} Coordinate-free isotropy is \(\mathsf R=\Sigma\mathsf g\propto I\). This turns the earlier numerical blindness of the commutator invariant into a two-line theorem: against isotropic noise the signal is exactly zero.

\bigskip
\textbf{Corollary 4.1.B (exponent invariance).} \emph{\(k(s)=\mathbb E_{\hat{\mathbf u}}\|B\hat{\mathbf u}\|^s\) contains no \(\Sigma\); twisting \(\Sigma\) by \(O\) cannot move the root of \(k(\cdot)=1\), so \(\alpha_\star(B^\top)=\alpha_\star(B)\).} This is a third route to the exponent invariance of §3, the noise-frame route. It is complementary rather than fully independent: the noise-frame route and the isospectral-Gram route both factor through \(k(s)\), while Proposition 3.5 is the independent statement.

\bigskip
\textbf{Corollary 4.1.C (directional profile, body frame).} \emph{Whenever the directional limits
\[
C^{\mathrm{body}}_{\mathbf n}(B;\Sigma):=
\lim_{u\to\infty}u^{\alpha_\star}
\mathbb P(\langle\mathbf y_t,\mathbf n\rangle>u)
\]
exist, the stationary body process obeys the polar-twist covariance
\[
\mathbf y^{B^\top,\Sigma}\stackrel d=O\,\mathbf y^{B,O^\top\Sigma O},
\qquad
C^{\mathrm{body}}_{\mathbf n}(B^\top;\Sigma)
=C^{\mathrm{body}}_{O^\top\mathbf n}(B;O^\top\Sigma O).
\]
The lab stationary law is isotropic under fresh Haar, so the nontrivial directional profile lives in the body frame.} The Haar-driver relabelling and complete proof are in Supplement S.B.3.

\bigskip
\textbf{Corollary 4.1.D (chiral null locus).} \emph{The commuting-twist set \(\{(B,\Sigma):[O(B),\Sigma]=0\}\) is a nontrivial null locus containing strata strictly larger than isotropy, and it is falsifiable: an engineered commuting-twist pair must show \(\Delta\log C=0\).}

Simulations check Corollaries A and D as finite-sample nulls at the working points (isotropy null \(|z|=0.75\), twist-equivalence null \(|z|=0.68\); C15--C16). Corollary C is verified as an exact directional identity (C30). The pathwise mapping of Theorem 4.1 holds to machine precision (\(\|B^\top-OBO\|\le1.4\times10^{-15}\), C3; the mapped chain agrees to \(1.6\times10^{-14}\), C4).

\paragraph{4.1.2 The generator, and the sign-chamber structure}\label{the-generator-and-the-sign-chamber-structure}

For a computable first-order law, we need the generator of \(O\) as \(A\) turns on. In the eigenbasis of \(S=\operatorname{diag}(s_1,s_2,s_3)\), write \(H_{ij}=1/(s_i+s_j)\) for the Hadamard kernel and \(J:=\operatorname{sgn}(S)\).

\bigskip
\textbf{Lemma 4.2 (twist generator; positive-definite and signed forms).} \emph{For \(A\) skew and \(\kappa\to0\):} \emph{(i) if \(S\succ0\) (so \(O(S)=I\)), \(\log O(S+\kappa A)=-2\kappa\,(A\circ H)+O(\kappa^3)\), odd in \(A\);} \emph{(ii) in general, with \(|s_i|\) distinct and nonzero,} \[
O(S+\kappa A)=J\cdot\exp(-\kappa\widetilde W)+O(\kappa^2),\qquad
\widetilde W_{ij}=(1+J_iJ_j)\,(A\circ H)_{ij}.
\]

\emph{Proof spine.} Linearising the polar-factor equation on the SPD chart gives \(\{S,-\partial_\kappa O|_0\}=2A\), hence \(\partial_\kappa O|_0=-2(A\circ H)\); parity removes the quadratic term. Decomposition under \(X\mapsto JXJ\) gives the signed chamber factor. Full details are in Supplement S.B.8.

The signed form (ii) controls the indefinite regime. The \textbf{chamber factor} \((1+J_iJ_j)\) is \(2\) on \emph{same-sign} eigen-pairs and \(0\) on \emph{mixed-sign} pairs. One asymptotic and one exact consequence follow.

\bigskip
\textbf{Assumption R (\(C^1\)-regularity).} \emph{The map \((B,\Sigma)\mapsto\log C(B;\Sigma)\) is \(C^1\) on an open full-matrix neighbourhood of the working point in \(GL(3)\times\mathrm{Sym}^+(3)\), so that \(\nabla_B\log C\) exists in the skew directions used below.} This is a named hypothesis, not a consequence of existence or positivity of \(C\); its scope and sensitivity audit are in Supplement S.B.8.

\bigskip
\textbf{Theorem 4.3 (chamber decomposition).} \emph{Under the hypotheses of Lemma 4.2(ii): (a) mixed-sign pairs decouple from the twist generator at first order --- and, within the nonzero-eigenvalue chamber, only such pairs can have \(s_i+s_j\to0\) by cancellation, so that would-be resonance multiplies a structurally vanishing channel. (b) Under Assumption R, the first-order mirror response splits into two channels,} \[
\Delta'(0)=\underbrace{\big\langle A-JAJ,\ \nabla_B\log C\big\rangle}_{\text{mixed-sign: direct}}+\underbrace{\big\langle \widetilde W,\ [\Sigma,D]\big\rangle}_{\text{same-sign: twist}},\qquad D:=\partial_\Sigma\log C|_{A=0},
\] \emph{a mixed-sign} direct \emph{channel and a same-sign} twist \emph{channel. (c) For \(J=I\) and \(S\succ0\) this reduces to Lemma 4.2(i)/Theorem 4.5 (the mixed sector is empty). Algebraically, \(\mathbb L_S=\{S,\cdot\}\) on the skew sector has eigenvalues \(s_i+s_j\); at a rank-\((d-1)\) positive-semidefinite boundary these pair sums can remain nonzero, but this observation does not extend the invertible polar chart or the response theorem to singular \(B\).}

\emph{Algebraic reading.} The skew-sector chamber split is the \(\pm1\)-eigenspace decomposition of \(X\mapsto JXJ\): the \(+1\) (same-sign) sector carries the twist and the \(-1\) (mixed-sign) sector the direct response.

\emph{Proof spine.} Expand the polar twist about \(J\Sigma J\) and use \(C(X;J\Sigma J)=C(JXJ;\Sigma)\). The two eigenspaces of \(X\mapsto JXJ\) give the direct mixed-sign and same-sign twist terms. Supplement S.B.8 gives the full derivation.

For mixed-sign perturbations, the first-order statement extends to all orders.

\bigskip
\textbf{Theorem 4.4 (exact mixed-sector rigidity).} \emph{If \(\{J,A\}=0\) --- equivalently, \(A\) is supported on mixed-sign pairs --- then \(JA\) is symmetric, so \(J(S+A)=|S|+JA\) is symmetric, and whenever \(|S|+JA\succ0\) the polar factor is \(O(S+A)=J\) exactly, to all orders. Consequently \(C\big((S+A)^\top;\Sigma\big)=C(S+A;\,J\Sigma J)\) exactly: the twist kernel never enters the mixed sector, at any order, even when \(s_i+s_j\) is arbitrarily close to zero.}

\emph{Proof.} \((JA)^\top=A^\top J=-AJ=JA\) by \(\{J,A\}=0\), so \(JA\) is symmetric; then \(B=J\cdot(JB)\) with \(JB=|S|+JA\) symmetric positive definite, and the orthogonal polar factor of such a product is \(J\) by uniqueness of the polar decomposition. \(\blacksquare\)

This resolves the ``resonance'' question. The apparent singularity of the naive first-order formula at \(s_i+s_j=0\) is neither a physical divergence nor another mechanism; it is \emph{channel re-routing}. Within the nonzero-eigenvalue chamber, cancellation makes a pair sum zero only for a mixed-sign pair. There the twist channel is closed (\(\widetilde W=0\)), and the regular direct channel carries the chiral response. Same-sign denominators vanish only at a singular boundary, outside the invertible polar chart. The pole multiplies a structural zero.

\paragraph{4.1.3--4.1.4 Regularity and the resolved near-miss}\label{regularity-through-near-pairing}
The apparent pole is a channel artefact: its coefficient vanishes on the mixed-sign sector, where the regular direct term carries the response. A near-pairing sweep confirms regular response and rejects the extended \(1/(s_i+s_j)\) scaling. It also records the misleading single-point agreement behind that false hypothesis. Supplement S.B.4 gives the full sweep and near-miss audit (C28--C29; checkpoint index S.E.2). No polar-twist pair-sum resonance remains; a closed direct-channel gradient theory is open.

\begin{center}\rule{0.5\linewidth}{0.5pt}\end{center}

\subsubsection{4.2 The first-order response: Ward identity, exact confinement, and fluctuation--dissipation}\label{the-first-order-response-ward-identity-exact-confinement-and-fluctuationdissipation}

\paragraph{4.2.1 The chiral response as a Lyapunov Green's function}\label{the-chiral-response-as-a-lyapunov-greens-function}

At \(B=S\succ0\) with distinct spectrum, Lemma 4.2(i)'s twist generator turns Theorem 4.1 into a closed first-order law. Let \(\Delta(t)=\log C(S+tA;\Sigma)-\log C(S-tA;\Sigma)\) and \(G:=\nabla_A\Delta|_0=2\,\mathrm{Skew}(\nabla_B\log C)|_{B=S}\). This definition is free of any polar gauge.

\bigskip
\textbf{Theorem 4.5 (first-order chiral response).} \emph{Under Assumption R, \(\Delta\log C=\langle A,G\rangle+o(\|A\|)\) --- with the remainder sharpening to \(O(\|A\|^3)\) under \(C^3\), since mirror-oddness removes the quadratic term --- where} \[
G=2\,H\circ[\Sigma,D],\qquad D=\partial_\Sigma\log C\big|_{A=0},
\] \emph{where \(D\) is the symmetric noise-gradient of the tail constant, measurable by achiral \(\Sigma\)-bumps.}

The formula has a Green's-function form that explains its domain and links it to the chamber structure.

\bigskip
\textbf{Proposition 4.6 (chiral Ward identity; resolvent form).} \emph{Under Assumption R, at \(B=S\succ0\), \(\{S,\ \mathrm{Skew}\,\nabla_B\log C\}=[\Sigma,D]\), equivalently} \[
\Delta'(0)=2\big\langle \mathbb L_S^{-1}A,\ [\Sigma,D]\big\rangle,\qquad \mathbb L_S=\{S,\cdot\}\big|_{\mathrm{skew}},\ \text{eigenvalues }s_i+s_j.
\]

\emph{Proof spines and scope.} From Theorem 4.1 and Lemma 4.2(i), \(O=I-2t(A\circ H)+o(t)\), so \(\partial_t(O^\top\Sigma O)|_0=2[A\circ H,\Sigma]\); differentiation and the commutator--trace identity give Theorem 4.5. Applying the self-adjoint Lyapunov map and \(\{S,A\circ H\}=A\) gives Proposition 4.6. The full derivations, including the distinction between this dynamical skew constraint and the kinematic equivariance identity, are in Supplement S.B.8. On the SPD chart every \(s_i+s_j>0\); algebraic invertibility more generally requires only \(s_i+s_j\ne0\), but that observation alone does not extend the response law into an indefinite chamber. Theorems 4.3--4.4 show why the observable remains finite at the apparent pole.

\paragraph{4.2.2 Parameter-free cross-validation at canonical scale}\label{parameter-free-cross-validation-at-canonical-scale}

\begin{figure}[H]
\centering
\includegraphics[width=0.95\textwidth]{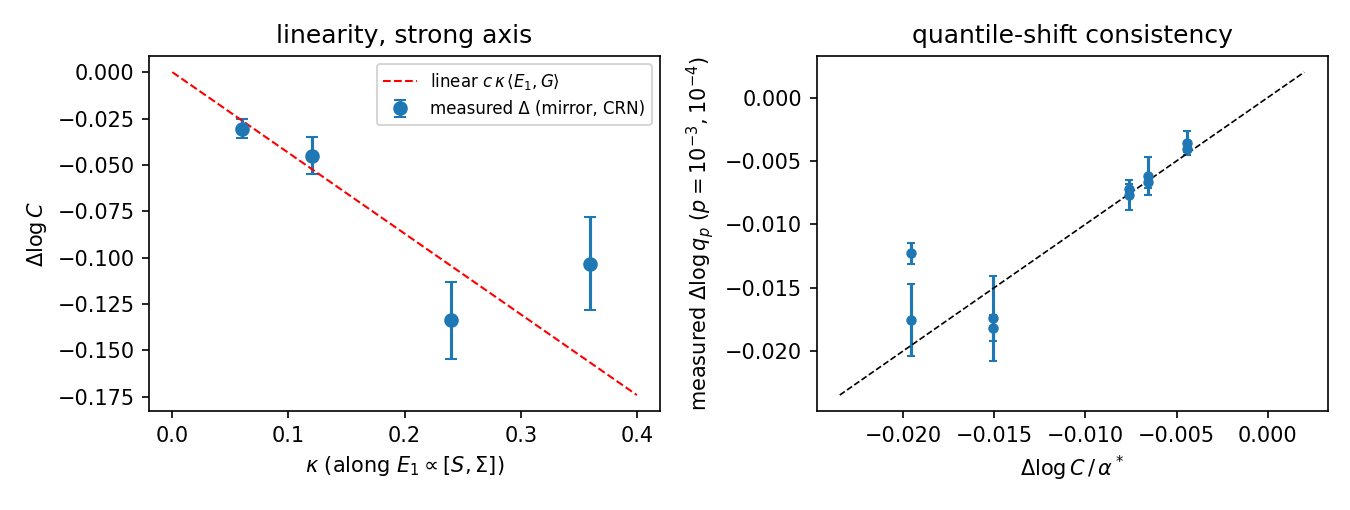}
\caption{First-order response and its level character (§4.2.2). \emph{Left:} in this \(\alpha_\star\)-recalibrated sweep, the mirror gap \(\Delta\log C\) is consistent with the predicted first-order line along the unit skew direction \(E_1=[S,\Sigma]/\|[S,\Sigma]\|_F\) of (4.4b) (dashed: \(\kappa\langle E_1,G_{\mathrm{win,pred}}\rangle\), no fitted constant); no departure is resolved through \(\kappa=0.30\). \emph{Right:} the measured log-quantile shift \(\Delta\log q_p\) (\(p=10^{-3},10^{-4}\); the mirror contrast \(dq_p\) of (2.6)) tracks \(\Delta\log C/\alpha_\star\) on the diagonal, consistent with a horizontal shift of a tail of fixed slope --- a level, not an exponent, effect.}\label{fig:f3-response}
\end{figure}

Figure~\ref{fig:f3-response} gives the main validation at a positive-definite point unlike the canonical one in \S2.4, calibrated to a \emph{near-canonical} index. Write
\[
S^+=d\cdot\operatorname{diag}(0.5,1.5,2.5),\qquad
\Sigma=R\,\operatorname{diag}(1.45,0.95,0.60)\,R^\top,\qquad
R=R_z(35^\circ)R_y(25^\circ)R_x(15^\circ),
\tag{4.4a}
\]
with \(d=0.516670153\) fixed by \(k(\alpha_\star)=1\) at \(\alpha_\star=6.836303\), giving \(s=(0.2583,0.7750,1.2917)\) and \(\|[S^+,\Sigma]\|_F=0.474\ne0\). Nine digits of \(d\) matter: \(\partial\alpha_\star/\partial d=-112\), so rounding to \(d=0.51667\) gives \(\alpha_\star=6.836320\). The index is near-canonical, not equal to the canonical operator's \(\alpha_\star=6.836931\).

The chiral sweep runs along the unit commutator direction
\[
E_1=\frac{[S^+,\Sigma]}{\|[S^+,\Sigma]\|_F},\qquad \|E_1\|_F=1,
\qquad
B^+(\kappa)=\lambda(\kappa)\bigl(S^++\kappa E_1\bigr),
\tag{4.4b}
\]
where \(\lambda(\kappa)\) is fixed at each \(\kappa\) by \(k(\alpha_\star)=1\). By homogeneity, \(\lambda(\kappa)=k(S^++\kappa E_1;\alpha_\star)^{-1/\alpha_\star}\). Thus \(\alpha_\star(\kappa)\equiv6.836303\), with \(\lambda\) running from \(0.99991\) at \(\kappa=0.02\) to \(0.98113\) at \(\kappa=0.30\). The quantile-to-level conversion uses this fixed exponent. Without rescaling, the index would fall to \(5.8195\) at \(\kappa=0.30\). Because it multiplies both parts, rescaling leaves \(\|A\|/\|S\|\) unchanged.

By Corollary 2.4, this point differs from \S2.4. The axial vector of \(E_1\), \(\boldsymbol\omega^+=(-0.1271,0.6643,0.2062)\), has three nonzero components in the simple eigenframe of \(S^+\), and
\[
\det[\boldsymbol\omega^+,S^+\boldsymbol\omega^+,(S^+)^2\boldsymbol\omega^+]=-4.80\times10^{-3}\ne0.
\]
Thus no \(P^+\in O(3)\) fixes \(S^+\) while reversing \(E_1\). The response operator is chiral in isolation, whereas the canonical operator satisfies Corollary 2.4(iii). Under Assumption R, the asymptotic objects are \(D_\infty:=\partial_\Sigma\log C|_{A=0}\) and \(G_\infty:=2H\circ[\Sigma,D_\infty]\). The controlled bumps instead determine a finite-window proxy \(D_{\mathrm{win}}\) and the parameter-free prediction \(G_{\mathrm{win,pred}}:=2H\circ[\Sigma,D_{\mathrm{win}}]\); Supplement S.B.5 gives the full design. Against the finite-window mirror-gradient proxy \(g_{\mathrm{win}}\),
\[
\cos\angle(g_{\mathrm{win}},G_{\mathrm{win,pred}})=0.9993\pm0.0019,
\qquad
\|g_{\mathrm{win}}\|/\|G_{\mathrm{win,pred}}\|=1.010\pm0.078.
\]
The sign agrees, and no departure from linearity in \(\kappa\) is resolved through \(\kappa=0.30\). The two tail windows differ by only \(0.004\pm0.010\), supporting a level shift at fixed slope. Supplement S.E.2 indexes the hard-null, leakage, and working-point audits.

\paragraph{4.2.3 Exact confinement and the leakage budget}\label{exact-confinement-and-the-leakage-budget}

The one-plane picture --- confinement to \(\{[S,\Sigma],[S^2,\Sigma]\}\) --- is a truncation. In \(d=3\), the exact result is a closed decomposition.

\bigskip
\textbf{Theorem 4.7 (exact confinement kernel).} \emph{In \(d=3\), under Assumption R, let \(S=S^\top\succ0\) have simple spectrum. Write \(D=q(S)+D_{\mathrm{off}}\), where \(q(t)=a_0+a_1t+a_2t^2\) is the unique polynomial of degree at most two satisfying \(q(s_i)=D_{ii}\). Then, exactly,} \[
G=-2a_2\,[S,\Sigma]\;-\;2a_1\,\mathbb L_S^{-1}[S,\Sigma]\;+\;2H\circ[\Sigma,D_{\mathrm{off}}],\qquad \mathbb L_S^{-1}[S,\Sigma]=H\circ[S,\Sigma],
\] \emph{where \(H_{ij}=(s_i+s_j)^{-1}\) off the diagonal. Equivalently, the diagonal-response kernel is \(\psi(s_i,s_j)=-2\big(a_2+a_1/(s_i+s_j)\big)\). The split is exact, not a Taylor truncation. A complete algebraic proof and component ledger are given in Supplement S.B.6.}

Under Assumption R, \(I_1\) overshoots the predicted first-order tail-level gradient at the validation point. At the tested threshold, it overshoots the measured finite-window proxy by \(37\%\). The resolvent term provides the required renormalisation. Supplement S.B.6 gives the exact component coordinates and leakage convention; S.E.2 indexes checkpoint C20. The finite closure is specific to \(d=3\); higher dimensions add commutator channels.

\paragraph{4.2.4 A fluctuation--dissipation representation of the noise-gradient}\label{a-fluctuationdissipation-representation-of-the-noise-gradient}

The noise gradient \(D=\partial_\Sigma\log C\) is a \emph{susceptibility}: a tail-conditioned innovation correlation. The order of limits matters because no fixed lag window carries the limiting response.

Write \(s(\boldsymbol\eta)=\tfrac12(\Sigma^{-1}\boldsymbol\eta\boldsymbol\eta^\top\Sigma^{-1}-\Sigma^{-1})\) for the one-step Gaussian score, and set \[
K(z):=\sum_{k\ge1}\Big(\mathbb E\big[\boldsymbol\eta_{-k}\boldsymbol\eta_{-k}^\top\mid r_0>z\big]-\Sigma\Big),\qquad
D(z):=\tfrac12\Sigma^{-1}K(z)\Sigma^{-1},
\] so that \(D(z)\) is the finite-threshold estimator of the susceptibility used in the controlled working-point diagnostics below.

\bigskip
\textbf{Assumption UI\(^{\prime\prime}\) (for Theorem 4.8).} \emph{(i) The untilted norm chain is \(V\)-geometrically ergodic and, for every fixed \(z\), the likelihood-ratio identity \(\partial_\Sigma\log\mathbb P(r_0>z)=\sum_{k\ge1}\mathbb E[s(\boldsymbol\eta_{-k})\mid r_0>z]\) holds with an absolutely convergent infinite-lag sum. (ii) As \(z\to\infty\), \(\partial_\Sigma\log\mathbb P(r_0>z)\to\partial_\Sigma\log C\) locally uniformly in \(\Sigma\) on compact subsets of \(\mathrm{Sym}^+(3)\).} Part (ii) is the substantive derivative/tail-limit exchange; it is separate from Assumption R and stronger than the \(\Sigma\)-differentiability component used there; it does not follow from \(\partial_\Sigma\alpha=0\). Scope and justification are audited in Supplement S.B.8; the coupling proof below and Supplement S.C.1 do not prove part (ii).

\bigskip
\textbf{Theorem 4.8 (tail-level fluctuation--dissipation; under Assumption UI″).} \emph{For the fresh-Haar norm chain with \(\boldsymbol\eta_t\sim\mathcal N(0,\Sigma)\),} \[
D=\partial_\Sigma\log C=\lim_{z\to\infty}\ \sum_{k\ge1}\mathbb E\big[s(\boldsymbol\eta_{-k})\,\big|\,r_0>z\big]=\lim_{z\to\infty}\tfrac12\Sigma^{-1}K(z)\Sigma^{-1}.
\] \emph{The limit must be taken outside the sum: for each fixed \(k\) the conditional score limit vanishes (Lemma 4.8.1(a)), so the interchanged object \(\sum_k\lim_z\) is identically zero and does not represent \(D\).}

Here \(\partial_\Sigma\alpha=0\) (§3) removes the divergent \((\partial_\Sigma\alpha)\log z\) term, allowing the tail limit to exist. Passing \(\partial_\Sigma\) \emph{through} the limit is Assumption UI″(ii). The exponent theorem thus supplies the second leg of the §3↔§4 bridge. Two results make this delocalisation precise \emph{without} Assumption UI″, under the standing model and Theorem LR's hypotheses.

\bigskip
\textbf{Lemma 4.8.1 (fixed-lag nullity and a uniform comparison rate).}
{\itshape Let
\[
\tau_k(z):=\mathbb E\!\left[
s(\boldsymbol\eta_{-k})\mid r_0>z
\right],
\qquad b:=\sigma_{\max}>1.
\]
(a) For every fixed \(k\ge1\),
\(\lim_{z\to\infty}\tau_k(z)=0\).

(b) Let \(\delta_0,C_\sharp,z_{\rm LR}\) be the constants in
Theorem LR (Appendix C; complete proof in Supplement S.C.2). Fix \(c_0>0\) with
\(c_0^2>2(\alpha_\star+1)\), put
\[
D_k(z):=
2c_0\sqrt{\lambda_{\max}(\Sigma)}\,
b^{k-1}\sqrt{\log z},
\]
and use the Frobenius norm on \(\mathrm{Sym}(3)\).
There are \(z_1<\infty\) and \(C_1<\infty\), depending only on
\((B,\Sigma,\alpha_\star,C,\varepsilon(\cdot),c_0)\), such that
\[
\|\tau_k(z)\|_F
\le
C_1 b^{k-1}\frac{(\log z)^{3/2}}{z}
\tag{4.8.1b}
\]
for every \(z\ge z_1\) and every \(k\ge1\) satisfying
\[
D_k(z)\le\delta_0z.
\tag{4.8.1c}
\]
In particular, for every fixed \(\varepsilon\in(0,1)\) there is
\(z_\varepsilon\) such that (4.8.1b) holds, with the same \(C_1\),
simultaneously for
\[
1\le k\le
(1-\varepsilon)\frac{\log z}{\log b},
\qquad z\ge z_\varepsilon.
\]}

\emph{Proof.}
For (a), condition on the most recent \(k\) drivers. Their \(k\)-step radial map satisfies
\[
|F_k(r)-P_kr|\le H_k,\qquad
P_k=\prod_{j=1}^k\|B\hat{\mathbf u}_{-j}\|,
\]
with \(H_k<\infty\) for fixed \(k\). Regular variation and fixed-\(k\) Gaussian-moment domination give
\[
\frac{\mathbb E[s(\boldsymbol\eta_{-k})\mathbf1_{\{r_0>z\}}]}{\bar F(z)}
\longrightarrow
\mathbb E[s(\boldsymbol\eta_{-k})]\,\mathbb E[P_k^{\alpha_\star}]=0.
\]

For (b), let \(\boldsymbol\eta'\) be an independent copy of the scored innovation and construct \(r'_0\) by replacing only \(\boldsymbol\eta_{-k}\) with \(\boldsymbol\eta'\), keeping all other drivers fixed. Then \(r'_0\stackrel d=R\) and is independent of the original scored innovation, hence
\[
N_k(z):=\mathbb E[s(\boldsymbol\eta_{-k})\mathbf1_{\{r_0>z\}}]
=\mathbb E\!\left[s(\boldsymbol\eta_{-k})
  \bigl(\mathbf1_{\{r_0>z\}}-\mathbf1_{\{r'_0>z\}}\bigr)\right].
\]
The one-step radial map is \(\|B\hat{\mathbf u}\|\)-Lipschitz, so the single replacement propagates as
\[
|r_0-r'_0|\le b^{k-1}\|\boldsymbol\eta-\boldsymbol\eta'\|.
\]
On the Gaussian good event
\[
G_z=\{\|\Sigma^{-1/2}\boldsymbol\eta\|^2,\,
       \|\Sigma^{-1/2}\boldsymbol\eta'\|^2\le c_0^2\log z\},
\]
the score is \(O(\log z)\) and the displacement is at most \(D_k(z)\). Indicator disagreement therefore implies \(|r'_0-z|\le D_k(z)\); under (4.8.1c), Theorem LR bounds this shell by
\[
C_\sharp D_k(z)\frac{\bar F(z)}z.
\]
The score-weighted complement of \(G_z\) is
\(O((\log z)^{3/2}z^{-c_0^2/2})\), uniformly in \(k\). Since
\(c_0^2/2>\alpha_\star+1\) and \(\bar F(z)\sim Cz^{-\alpha_\star}\), division by \(\bar F(z)\) yields (4.8.1b), with a constant independent of \(k,z\). Finally,
\[
k\le(1-\varepsilon)\frac{\log z}{\log b}
\quad\Longrightarrow\quad
\frac{D_k(z)}z=O\!\left(z^{-\varepsilon}\sqrt{\log z}\right)\to0,
\]
which gives the displayed uniform window. The explicit constants, \(\chi^2_3\) tail integration, and the unused \(H_k\) route are audited in Supplement S.C.1. \(\square\)

\bigskip
\textbf{Corollary 4.8.2 (lag-window escape).}
{\itshape Let \(K(z)\) be integer-valued.  If
\[
\frac{b^{K(z)}(\log z)^{3/2}}z\longrightarrow0,
\tag{4.8.2a}
\]
then
\[
\sum_{k=1}^{K(z)}
\|\tau_k(z)\|_F\longrightarrow0.
\tag{4.8.2b}
\]
In particular, for every fixed \(\varepsilon\in(0,1)\),
(4.8.2b) holds whenever
\[
K(z)\le
(1-\varepsilon)\frac{\log z}{\log b}.
\]
More sharply, it holds for
\[
K(z)=
\left\lfloor
\frac{
\log z-\frac32\log\log z-\omega(z)}
{\log b}
\right\rfloor
\quad\text{for every }\omega(z)\to+\infty,
\]
provided the numerator tends to \(+\infty\) if a nonempty growing
window is intended.}

\emph{Proof.}
The shell-domain condition of Lemma 4.8.1(b) is automatic under
(4.8.2a), because
\[
\frac{b^{K(z)}\sqrt{\log z}}z
=
\frac1{\log z}\,
\frac{b^{K(z)}(\log z)^{3/2}}z
\longrightarrow0.
\]
Therefore
\[
\begin{aligned}
\sum_{k=1}^{K(z)}\|\tau_k(z)\|_F
&\le
C_1\frac{(\log z)^{3/2}}z
\sum_{k=1}^{K(z)}b^{k-1}\\
&=
\frac{C_1}{b-1}
\frac{b^{K(z)}-1}{z}
(\log z)^{3/2}\\
&\le
\frac{C_1}{b-1}
\frac{b^{K(z)}(\log z)^{3/2}}z
\longrightarrow0.
\end{aligned}
\]
For \(K(z)\le(1-\varepsilon)\log z/\log b\), the last display is
\(O(z^{-\varepsilon}(\log z)^{3/2})\). \(\square\)

The susceptibility escapes every window \(K(z)\) satisfying \(b^{K(z)}(\log z)^{3/2}/z\to0\), including each fixed-fraction window \(K(z)\le(1-\varepsilon)\log z/\log b\) and the displayed log--log refinement. Its contributing lags therefore recede as the tail deepens. Its natural age scale is heuristically \(k\sim\log z/\Lambda'(\alpha_\star)\), consistent with the proven floor because \(\Lambda'(\alpha_\star)\le\log\sigma_{\max}\); a closed-form limiting lag profile is Open Problem 1, §4.7. The delocalisation also has an exact finite-\(z\) constraint, checking the representation and controlled working-point score estimator.

\bigskip
\textbf{Proposition 4.8.3 (finite-\(z\) sum rule; under Assumption UI″(i)).} \emph{With \(\alpha_{\mathrm{loc}}(z):=-z\,\partial_z\log\mathbb P(r_0>z)\) the local tail slope, \(2\langle D(z),\Sigma\rangle=\langle K(z),\Sigma^{-1}\rangle=\alpha_{\mathrm{loc}}(z)\), and \(\alpha_{\mathrm{loc}}(z)\to\alpha_\star\) under Assumption UI″(ii).}

\emph{Proof.} The survival law scales exactly, \(\mathbb P_{\lambda\Sigma}(r_0>z)=\mathbb P_\Sigma(r_0>z\lambda^{-1/2})\); differentiating at \(\lambda=1\), \(\partial_{\log\lambda}\log\mathbb P_{\lambda\Sigma}(r_0>z)|_{\lambda=1}=+\tfrac12\alpha_{\mathrm{loc}}(z)\) --- positive, since increasing the innovation scale raises a fixed-threshold exceedance probability. By UI″(i) the left side is \(\langle D(z),\Sigma\rangle\); the limit follows from UI″(ii) and \(C(\lambda\Sigma)=\lambda^{\alpha_\star/2}C(\Sigma)\). \(\qquad\square\)

In the controlled working-point simulation, \(2\langle D(z),\Sigma\rangle\) is \(6.13\) at the shallow threshold and \(6.49\) at the deeper one, approaching \(\alpha_\star=6.836\). This independently checks any closed form for the limiting lag profile.

\emph{Sketch of Theorem 4.8.} The finite-\(z\) step is the classical likelihood-ratio / score-function gradient (Glynn 1990; L'Ecuyer 1990): under UI″(i), \(\partial_\Sigma\log\mathbb P(r_0>z)=\sum_k\mathbb E[s(\boldsymbol\eta_{-k})\mid r_0>z]\), an absolutely convergent infinite-lag sum at each fixed \(z\). Writing \(\log C=\lim_z(\log\mathbb P(r_0>z)+\alpha\log z)\) and using \(\partial_\Sigma\alpha\equiv0\) (§3) gives \(D=\lim_z\partial_\Sigma\log\mathbb P(r_0>z)\), the limit outside the sum; the interchange with \(\partial_\Sigma\) is UI″(ii). That the mass does not stay at the excursion seed is Lemma 4.8.1 and Corollary 4.8.2 --- each fixed-lag term vanishes and the response escapes every finite window, which is exactly why the naive \(\sum_k\lim_z\) is zero and the correct object needs the outer limit. \(\qquad\square\)

\bigskip
\textbf{Working-point heuristic (sign-blind seed term).} \emph{The \(\alpha\)-tilted direction law is explicit, \(\nu_\alpha(d\hat{\mathbf u})\propto\|S\hat{\mathbf u}\|^\alpha\sigma(d\hat{\mathbf u})\), and the seed-noise contribution concentrates on top-\(|S|\) directions. This motivates, but does not prove, a \(|s_i|\)-ordering for the total susceptibility \(D\): the unclosed stationary-law response from \(k\ge2\) can alter it.}

At two controlled thresholds, the conditional-score and covariance-bump routes agree in direction (\(\cos=0.9992,0.9976\)). Their magnitude ratio rises \(0.907\to0.967\), consistent with the documented local-slope conversion. At the indefinite working point, the observed order \(D_{33}>D_{11}>D_{22}\) follows \(|s_3|>|s_1|>|s_2|\), but is a numerical pattern, not a monotonicity theorem. Supplements S.C.1 and S.E.2 give the conversion and checkpoint details.

Two negative results sharpen the closed-form question. First, ``the explicit Goldie term equals the FDT lag-1 term'' is false: \(\|D_{\mathrm{lag1}}\|/\|D_{\mathrm{exp}}\|=0.066\text{–}0.083\) (a \(12\text{–}15\times\) mismatch) at \(\cos=0.885\). The Goldie split \(D=D_{\mathrm{exp}}+D_\pi\) --- the explicit final-innovation term \(D_{\mathrm{exp}}\) plus stationary-law response \(D_\pi\) --- and the FDT lag sum are \emph{inequivalent decompositions of the same total}; only their totals must agree (Remark 4.12, C21). Second, the approximation \(D\approx D_{\mathrm{exp}}\) fails quantitatively (\(55\%\) of the norm, \(31^\circ\) in direction). Thus \(D_\pi\) in the Goldie route, or the \(k\ge2\) tilted sum in the FDT route, remains the hard part of any closed form for \(\varphi\) and is left open (§4.7).

\paragraph{4.2.5 The §3↔§4 bridge}\label{the-34-bridge}

At the canonical point for §3's exponent theory --- \(S=\operatorname{diag}(-c,0,c)\), diagonal \(\Sigma\) --- the level effect also vanishes by an exact symmetry requiring neither regularity nor positive-definiteness.

\bigskip
\textbf{Theorem 4.9 (mirror-evenness at commuting points).} \emph{Let \([S,\Sigma]=0\) (any signature, any rank --- the canonical \(S=\operatorname{diag}(-c,0,c)\) included), and in the shared eigenbasis write \(A=a_{12}E_{12}+a_{13}E_{13}+a_{23}E_{23}\). Then \(\Delta\log C(\kappa A)\) is odd in each coordinate \(a_{ij}\) separately. Consequently: (i) if any one coordinate vanishes, \(\Delta\log C(\kappa A)\equiv0\) for every \(\kappa\) for which the §3 tail constants exist, exactly and without Assumption R; (ii) if \(\Delta\log C\) is \(C^3\) at the point, then \(\Delta\log C(\kappa A)=c\,a_{12}a_{13}a_{23}\kappa^3+o(\kappa^3)\); if it is analytic, every nonzero monomial is odd in all three coordinates.}

\emph{Proof.} For the sign-flip \(Q_i=\operatorname{diag}(\ldots,-1,\ldots)\) (the \(-1\) in slot \(i\)), \(Q_iSQ_i=S\) and \(Q_i\Sigma Q_i=\Sigma\) (both diagonal), so conjugation by \(Q_i\) is a symmetry of the noise law, giving \(\Delta\log C(Q_iAQ_i)=\Delta\log C(A)\). But \(Q_iAQ_i\) flips the sign of exactly the two coordinates carrying index \(i\); combined with mirror-oddness \(\Delta\log C(-A)=-\Delta\log C(A)\), this forces \(\Delta\log C\) to be odd in each coordinate separately. (i) and (ii) follow. \(\blacksquare\)

Part (i) needs neither differentiability nor a polar chart, but applies only where both §3 tail constants exist. Numerically, when \(A\) has a vanishing pair-coordinate, the CRN-paired chains for \(S\pm\kappa A\) coincide bitwise (C31, tested for the \((a,0,c)\) and \((0,b,c)\) families). Under the \(C^3\) condition, a generic response has cubic leading Taylor order. Part (ii) \emph{predicts} the cubic residual at commuting working points: the leading nonzero response is trilinear in the three chiral coordinates, hence \(O(\kappa^3)\). On the positive-definite side, the first-order formula gives the same content --- \([S,\Sigma]=0\Rightarrow D\) diagonal in the shared basis (by the same sign-flip stabiliser) \(\Rightarrow[\Sigma,D]=0\Rightarrow G=0\) --- but Theorem 4.9 also holds outside Theorem 4.5's chart.

The two legs \emph{meet} on the commuting locus for different structural reasons. A general theorem protects the \textbf{exponent} (Proposition 3.5 --- every i.i.d. ensemble, every direction, in the \(\Lambda\)-root regime). For the \textbf{level}, the displayed alignments give sufficient nulls; they do not exclude other or accidental zeros. At the tested rotated point, the finite-threshold contrast is consistent with a nonzero linear coefficient. Identifying it with the asymptotic derivative remains conditional on Assumption R. At the canonical point, both displayed channels are exactly hidden; §3 and §4 give the same statement from two directions.

\begin{center}\rule{0.5\linewidth}{0.5pt}\end{center}

\subsubsection{4.3 What the lab frame hides: one-lag operator collapse, frame-invariance of the mirror gap, and the fourth-moment channel}\label{what-the-lab-frame-hides-one-lag-operator-collapse-frame-invariance-of-the-mirror-gap-and-the-fourth-moment-channel}

Section 4.1 has a precise observability consequence. Switching between laboratory and body representations of the same strict-Haar model does not change the possibly zero mirror gap. Strict-Haar averaging only prevents recovery of the transpose-odd operator from a \emph{one-lag linear} fit; the radial mirror gap remains unchanged. Three results establish this: one-lag linear collapse (Theorem 4.10), frame-invariance of the mirror gap (Corollary 4.10.A), and survival of transpose-odd content in the fourth radial moment (Proposition 4.10.B).

\bigskip
\textbf{Theorem 4.10 (strict-Haar collapse of the one-lag linear operator).} \emph{Under a strict-Haar rotation law, \(\mathbb E_{\mathrm{Haar}}[Q B Q^\top]=(\operatorname{tr}B/3)\,I\). Consequently the population one-lag regression coefficient of the strict-Haar laboratory process is \((\operatorname{tr}B/3)\,I\): its antisymmetric part vanishes identically. Any lab-frame statistic that depends only on the antisymmetric part of the one-lag} linear \emph{predictor therefore carries zero population signal --- the transpose-odd part of \(B\) is invisible to one-lag linear observation.} (The averaging identity is the standard Schur/Haar first moment; the content is its consequence for one-lag linear observability.)

\bigskip
\textbf{Corollary 4.10.D (all-lag linear collapse).} \emph{The same averaging fixes every lag: for \(k\ge1\),} \[
\mathbb E[\mathbf x_{t+k}\mathbf x_t^\top]=\Big(\tfrac{\operatorname{tr}B}{3}\Big)^{k}\,\tfrac{m_2}{3}\,I,\qquad m_2:=\mathbb E\|\mathbf x_t\|^2=\frac{\operatorname{tr}\Sigma}{1-\|B\|_F^2/3}
\] \emph{(finite whenever \(\|B\|_F^2<3\)). Every lag-\(k\) cross-moment is therefore isotropic, and the entire second-order (Gaussian) structure of the lab process is fixed by the two scalars \((\operatorname{tr}B,\,m_2)\) --- both transpose-even. When \(S\) is traceless the collapse is total: \(\operatorname{tr}B=\operatorname{tr}S=0\) forces every lag-\(k\) (\(k\ge1\)) cross-moment to vanish identically, so the linear channel is then empty rather than only transpose-blind.}

\emph{Proof.} The noise cross-terms vanish (\(\mathbb E[Q_j\boldsymbol\eta_j]=0\), independence of the fresh step); the operator factor is multiplicative across the \(k\) independent steps, \(\mathbb E[\prod_{j=1}^{k}Q_jBQ_j^\top]=(\mathbb E[QBQ^\top])^{k}=(\operatorname{tr}B/3)^{k}I\) by Theorem 4.10; and the stationary law is isotropic, \(\mathbb E[\mathbf x\mathbf x^\top]=(m_2/3)I\), with \(m_2\) fixed by the stationarity balance \(m_2/3=(\|B\|_F^2/3)(m_2/3)+\operatorname{tr}\Sigma/3\). \(\qquad\blacksquare\)

The linear channel lacks transpose-odd content at every lag, not only lag one. The Gaussian sector retains the isotropic autocovariance \((\operatorname{tr}B/3)^{k}(m_2/3)I\), which reduces to the lag-\(0\) variance \(m_2/3\) when \(S\) is traceless. It carries chirality magnitude only \emph{unidentifiably}, through \(\|B\|_F^2=\|S\|_F^2+\|A\|_F^2\), inseparable from \(\|S\|_F^2\) and \(\operatorname{tr}\Sigma\). No lag-\(k\) cross-moment --- equivalently, no population linear least-squares predictor --- distinguishes \(B\) from \(B^\top\). This conclusion concerns only \emph{second-order} structure and gives no distributional separator between the full laws under \(B\) and \(B^\top\). When both radial tail levels exist, their gap depends on the \emph{pair} of radial laws and may vanish or not. The one-sample sign symmetry \(\mathbf x_\infty\stackrel d=-\mathbf x_\infty\) is not a \(B\leftrightarrow B^\top\) comparison. A valid one-sample nonlinear or cross-tail contrast remains open.

This concerns the one-lag \emph{linear} estimator, not the stationary law. Two results fix its scope: the mirror gap is unchanged (Corollary 4.10.A), and transpose-odd content returns at the next radial moment (Proposition 4.10.B).

\bigskip
\textbf{Corollary 4.10.A (frame-invariance of the mirror gap).} \emph{The mirror gap \(\Delta\log C=\log C(B;\Sigma)-\log C(B^\top;\Sigma)\) is a functional of the pair of radial tail constants for \(B\) and \(B^\top\), and it may vanish (Theorem 3.3). Because the co-moving norm process is frame-free --- \(\|\mathbf x_t\|=\|\mathbf y_t\|\) pathwise (§2.1) and \(\|R\mathbf x\|=\|\mathbf x\|\) for every \(R\in SO(3)\) --- each radial law and tail constant, hence their pair and log-gap, is represented identically whether magnitudes are read in the laboratory frame \(\mathbf x_t\) or the co-moving body frame \(\mathbf y_t\). It is unchanged by a global \(\mathsf g\)-orthogonal change of observation frame, and is covariantly re-expressed under a general \(GL(3)\) basis change only when \(\mathsf g\) is carried along. The strict-Haar lab therefore does not dilute the mirror gap; it is the identical radial number in the lab and body descriptions of the same fixed model.}

\emph{Proof.} The tail constant is defined by \(\mathbb P(\|\mathbf x_\infty\|>u)\sim C\,u^{-\alpha_\star}\) (Theorem 3.3), a functional of the magnitude only; the radial reduction (Lemma 3.1) writes the magnitude chain as \(r_{t+1}=\|r_t B\mathbf V_t+\boldsymbol\eta_t\|\), which carries no frame label since \(\|\cdot\|\) is rotation-invariant. Concretely \(\|\mathbf x_t\|=\|\mathbf y_t\|\) for every realisation (§2.1), so the lab and co-moving body frames induce the identical radial law; hence \(C\), the pair \((C(B;\Sigma),C(B^\top;\Sigma))\), and their log-ratio \(\Delta\log C\) coincide across the two frames. \(\qquad\blacksquare\)

\emph{Scope.} This is invariance under a \textbf{change of observation frame}: a norm-preserving global rotation applied simultaneously to \((\mathbf x_t,B,\Sigma)\), together with equality of the lab and co-moving radial laws. It is \textbf{not} invariance under rotating \(\Sigma\) \emph{relative to} \(B\). The twist \(\Sigma\mapsto O^\top\Sigma O\), with \(B\) fixed, moves \(C\) (Theorem 4.1) and supplies the mechanism of §4. Two other changes lie outside this invariance. One changes the \emph{noise-transport rule} from co-rotating innovations to laboratory-fixed covariance; in the §2.1 control, the mirror contrast becomes a powered null rather than frame-invariant. The other changes the frame process's \emph{temporal law}; in a controlled refresh-rate family, \(\gamma\) and the effective index \(\hat\alpha\) vary with \(\omega\). Both change the model or observable. The result concerns only lab/body representations of one fixed model: representation invariance, not dynamical invariance.

\bigskip
\textbf{Proposition 4.10.B (the transpose-odd content has surviving channels at fourth radial moment).} \emph{Theorem 4.10 annihilates the transpose-odd part of the one-lag} linear \emph{predictor, but not of the one-lag law. The vector conditional mean \(\mathbb E[\mathbf x_{t+1}\mid\mathbf x_t]=(\operatorname{tr}B/3)\mathbf x_t\) and the radial second moment \(\mathbb E[\|\mathbf x_{t+1}\|^2\mid\|\mathbf x_t\|=r]\) are mirror-even (the latter because \(\mathbb E\|B\hat{\mathbf u}\|^2\) is a singular-value functional). Moments of orders one and three need not be mirror-even; among even integer orders, the fourth is the first whose mirror difference can be nonzero, with a closed coefficient:} \[
\Delta m_4(r):=\mathbb E\big[\|\mathbf x_{t+1}\|^4\,\big|\,\|\mathbf x_t\|=r;B\big]-\mathbb E\big[\cdots;B^\top\big]=-\tfrac{8r^2}{3}\,I_1,\qquad I_1=\langle A,[S,\Sigma]\rangle,
\] \emph{the commutator invariant of §4.2.3.}

\emph{Proof.} Since \(Q_t\) is orthogonal, \(\|\mathbf x_{t+1}\|=\|r B\hat{\mathbf u}+\boldsymbol\eta\|\) with \(r=\|\mathbf x_t\|\), \(\hat{\mathbf u}=Q_t^\top\mathbf x_t/r\) uniform on the sphere and \(\boldsymbol\eta\sim\mathcal N(0,\Sigma)\) independent. For fixed \(\mathbf w=rB\hat{\mathbf u}\), the Gaussian quartic identity gives \(\mathbb E_{\boldsymbol\eta}\|\mathbf w+\boldsymbol\eta\|^4=\|\mathbf w\|^4+2\|\mathbf w\|^2\operatorname{tr}\Sigma+(\operatorname{tr}\Sigma)^2+4\,\mathbf w^\top\Sigma\mathbf w+2\operatorname{tr}\Sigma^2\). Under \(B\mapsto B^\top\) the terms in \(\|\mathbf w\|^2=r^2\hat{\mathbf u}^\top B^\top B\hat{\mathbf u}\) and \(\|\mathbf w\|^4\) are unchanged (they depend on \(B\) only through \(B^\top B\), isospectral with \(B B^\top\) under the sphere average, using \(\mathbb E[\hat{\mathbf u}\hat{\mathbf u}^\top]=I/3\) and \(\mathbb E[(\hat{\mathbf u}^\top G\hat{\mathbf u})^2]=[(\operatorname{tr}G)^2+2\operatorname{tr}G^2]/15\)); the only transpose-odd term is the cross term \(4r^2\,\mathbb E_{\hat{\mathbf u}}[\hat{\mathbf u}^\top B^\top\Sigma B\hat{\mathbf u}]=\tfrac{4r^2}{3}\operatorname{tr}(B^\top\Sigma B)=\tfrac{4r^2}{3}\operatorname{tr}(\Sigma B B^\top)\), and \(BB^\top=S^2-A^2+[A,S]\) contributes the odd piece \(\tfrac{4r^2}{3}\operatorname{tr}(\Sigma[A,S])\). Subtracting the mirror (\([A,S]\mapsto-[A,S]\)) doubles it: \(\Delta m_4=\tfrac{8r^2}{3}\operatorname{tr}(\Sigma[A,S])=-\tfrac{8r^2}{3}\langle A,[S,\Sigma]\rangle\). No stationary fourth moment is required --- this is the one-step conditional law. \(\qquad\blacksquare\)

Theorem 4.10's blindness is specific to the \emph{linear} (vector-mean / radial-second-moment) channel; a \emph{radial-nonlinear} one-lag statistic retains transpose-odd content. Two conclusions follow. First, \(I_1\) vanishes when \([S,\Sigma]=0\), the \emph{commuting} locus. This locus includes but exceeds isotropic \(\Sigma\): a \(\Sigma\) diagonal in the eigenbasis of \(S\) commutes without being isotropic (cf.~Theorem 4.9's counterexample). This explains why the modeled one-step channel weakens as residual covariance approaches the commuting locus. Second, the \emph{first-order level} response renormalises \(I_1\) by the resolvent. Under Assumption R, \(I_1\) overshoots the predicted first-order tail-level gradient and, at the tested threshold, the measured proxy by \(37\%\) (§4.2.3). The \emph{fourth-moment one-step} response instead equals \(I_1\) without renormalisation. They are different observables of the same commutator.

Two cautions limit estimators based on this channel. At stationarity, \(\|\mathbf x_{t+1}\|^4\) has tail index \(\alpha_\star/4=1.709<2\) and infinite variance, ruling out ordinary finite-variance or \(\sqrt n\) inference. This differs from the signed Goldie functional in Lemma 4.13, for which only \(L^1\setminus L^2\) is established. The dependent chain needs separate sampling theory. The channel also has a blind direction: when \(A\perp[S,\Sigma]\) in the trace inner product, \(\Delta m_4\equiv0\) while \(\Delta\log C\) may remain nonzero. Thus a null statistic need not imply a null effect. A practical fourth-moment estimator and its finite-sample behavior under blind frame reconstruction remain open (§4.7).

A third caution confines this channel to a \emph{magnitude} reading on frame-blind data. Proposition 4.10.B concerns the mirror \emph{difference} \(\Delta m_4\), which subtracts \(B^\top\) and cancels all transpose-even terms. Data provide \(B\), not its mirror. The controlled contrast instead uses a same-\(S\) symmetric surrogate (\(A=0\), the same \(S,\Sigma\)), so it is not clean. Let \(c_2^{(k)}\) be the coefficient of \(r^2\) in \(\mathbb E[\|\mathbf x_{t+1}\|^{2k}\mid\|\mathbf x_t\|=r]\). The surrogate-differenced fourth-moment coefficient is
\[
c_2^{(2)}(S+A;\Sigma)-c_2^{(2)}(S;\Sigma)=\underbrace{-\tfrac43\langle A,[S,\Sigma]\rangle}_{\text{transpose-odd},\ O(\|A\|)}\;+\;\underbrace{\big[\tfrac23\operatorname{tr}\Sigma\,\|A\|_F^2-\tfrac43\operatorname{tr}(\Sigma A^2)\big]}_{\text{transpose-even},\ O(\|A\|^2)} .
\]
The odd part is exactly the mirror content of Proposition 4.10.B. The even part is a transpose-\emph{blind} \(\|A\|^2\) magnitude removed by the mirror difference but retained by the surrogate difference. At the working point, it dominates along the tested direction for \(\|A\|_F\gtrsim0.02\) (the ratio reaches \(\approx35\) at the working \(\kappa\)). Within the fresh-Haar spherical-input model, Proposition 4.10.C gives only an operator-uniform obstruction: no finite scalar-weighted radial combination cancels the even quadratic form for every operator direction while retaining the odd covector. This does not exclude cancellation tuned to one fixed operator-dependent direction or a jointly identified statistic using the operator.

\bigskip
\textbf{Proposition 4.10.C (operator-uniform radial sum-rule obstruction).} \emph{For an integer \(k\ge1\), assume \(\mathbb E\|\boldsymbol\eta\|^{2k-2}<\infty\), write \(A=[\boldsymbol\omega]_\times\), and let \(c_2^{(k)}(B)\) denote the expectation of the coefficient of \(r^2\) in the polynomial \(r\mapsto\|rB\hat{\mathbf u}+\boldsymbol\eta\|^{2k}\) (finite under this moment condition). For \(k\ge2\), set}
\[
T_k=\mathbb E[\boldsymbol\eta\boldsymbol\eta^\top\|\boldsymbol\eta\|^{2k-4}],\qquad c_k=\frac{2k(k-1)}3,\qquad a_k=\frac{2k^2}{3}\mathbb E\|\boldsymbol\eta\|^{2k-2},
\]
\emph{and define \(M_k=a_kI-c_kT_k\), \(v_k=-2c_k\operatorname{axl}[S,T_k]\). For \(k=1\), set \(M_1=(2/3)I\), \(v_1=0\). Then}
\[
c_2^{(k)}(S+A)-c_2^{(k)}(S)=\boldsymbol\omega^\top M_k\boldsymbol\omega+\boldsymbol\omega^\top v_k,
\qquad
\sum_k w_kM_k=0\ \Longrightarrow\ \sum_k w_kv_k=0.
\]
\emph{Thus no finite scalar-weighted combination can cancel the transpose-even quadratic form uniformly over operator directions while retaining the transpose-odd covector. Fixed-direction, operator-dependent cancellation is not excluded.}

\emph{Proof.} For \(k\ge2\), the coefficient of \(r^2\) in \((\|\eta\|^2+2r\eta^\top B\hat{\mathbf u}+r^2\|B\hat{\mathbf u}\|^2)^k\), followed by spherical averaging, is
\[
c_2^{(k)}(B)=\frac{k}{3}\mathbb E\|\eta\|^{2k-2}\|B\|_F^2+c_k\operatorname{tr}(BB^\top T_k).
\]
Using \(BB^\top=S^2-A^2+[A,S]\), \(A^2=\boldsymbol\omega\boldsymbol\omega^\top-\|\boldsymbol\omega\|^2I\), and \(\|A\|_F^2=2\|\boldsymbol\omega\|^2\) gives the stated \(M_k,v_k\). For \(k=1\), direct spherical averaging gives \(c_2^{(1)}(B)=\|B\|_F^2/3\). Finally, \(\sum_kw_kM_k=0\) implies \(\sum_{k\ge2}w_kc_kT_k=(\tfrac23w_1+\sum_{k\ge2}w_ka_k)I\); hence its commutator with \(S\), and therefore \(\sum_kw_kv_k\), vanishes. The expanded calculation and scope audit are in Supplement S.D.2. \(\qquad\blacksquare\)

No pseudoscalar escape exists: in \(d=3\), reflection is \(-1\) times a rotation, and centered Gaussian noise is sign-symmetric. Thus the strict-Haar lab law is \(O(3)\)-invariant and every parity-odd expectation vanishes. We do not use exploratory angular moment banks as evidence. The transpose-even second-moment excess \(\|A\|_F^2/3\) is an exact magnitude contrast only against a known same-\(S\) symmetric surrogate; one frame-blind sample does not identify it. Sign requires the co-moving frame (§4.5), its recovery, or more operator information.

\bigskip
\textbf{Remark (parity of the chiral observables).} The two orientation channels are both odd. Under Assumption R, \(\Delta\log C(\kappa A)=\kappa\langle A,G\rangle+o(\kappa)\), while \(I_1\) is exactly linear. The same-\(S\) surrogate-magnitude contrast is even and quadratic. Within the scalar radial combinations of Proposition 4.10.C, orientation cannot be separated uniformly from the even magnitude without operator information; the proposition does not rule out every nonlinear frame-blind statistic. In the controlled same-\(S\) comparison, magnitude is a second-order question that is stable under sign reversal, whereas orientation is first-order, transpose-odd, and sign-sensitive. The exponent/level dichotomy of §3--§4 is the same parity seen once more: the exponent is an even (transpose-blind) invariant, the level an odd (transpose-covariant) response.

Three distinctions define the interpretation. \textbf{(i)} The lab frame does \emph{not} dilute the mirror gap \(\Delta\log C\). By Corollary 4.10.A, it is the same frame-invariant radial number in lab and co-moving coordinates. \textbf{(ii)} The strict-Haar lab removes only the ability to \emph{recover} the transpose-odd part of \(B\) from a one-lag \emph{linear} fit. Haar averaging maps the one-lag operator to \((\operatorname{tr}B/3)I\) and erases its skew part (Theorem 4.10). \textbf{(iii)} This removal is confined to the linear channel. Transpose-odd content survives in the fourth radial moment (Proposition 4.10.B), so the obstruction concerns \emph{linear} one-lag recovery, not observation itself. In an operator-informed same-\(S\) comparison, a radial-nonlinear contrast can measure magnitude. One frame-blind sample identifies neither magnitude nor sign without more structure. Orientation requires the co-moving frame (§4.5), its recovery, or operator information. The odd operator is not identifiable from single-step \emph{linear} structure. As in §3.3, transpose-odd content excluded from the exponent and a linear one-lag operator can persist in the level and nonlinear radial structure.

Between a frozen frame and the strict-Haar endpoint, attenuation depends on the model and estimator. Controlled scrambling checks the implementation against the proved fast-mixing endpoint. It supplies no confidence interval for a retention frontier, calibrated map to \(\rho\), or certified recovery threshold. We use the strict-Haar result only as a counterfactual boundary.

\subsubsection{4.4 The observable and its contamination (methods)}\label{the-observable-and-its-contamination-methods}

Recovering the chiral level from a co-moving process requires body-frame reconstruction. The estimator controls two resulting contaminations.

Here \(C_0\) is the stationary lag-0 covariance and \(G_1\) the lag-1 cross-moment of observed increments; \(\hat C_0,\hat G_1\) are their sample estimates. The target \(A_{\mathrm{eff}}\) is the effective transpose-odd operator, and \(\hat w_{\mathrm{res}}\) is the estimated residual drift.

\bigskip
\textbf{Proposition 4.11 (population target of the fitted linear proxy).} \emph{For the fitted one-lag proxy, the naively antisymmetrised population cross-moment is not \(A_{\mathrm{eff}}\) but is contaminated at leading order by \(\tfrac12[S,C_0]\) (a \(39\%\) contamination at the working point). The corrected sample estimator \(\operatorname{skew}(\hat G_1\hat C_0^{-1})\) targets the population identity \(\operatorname{skew}(G_1C_0^{-1})=A_{\mathrm{eff}}\).}

\emph{Proof.} For the fitted one-lag proxy \(x_{t+1}=Bx_t+\eta_t\) with \(B=S+A\) and \(\mathbb E[\eta_t x_t^\top]=0\), the one-lag cross-moment is \(G_1=\mathbb E[x_{t+1}x_t^\top]=BC_0\) with \(C_0=\mathbb E[x_tx_t^\top]=C_0^\top\). Right-multiplying by \(C_0^{-1}\) returns \(B\) exactly, so \(\operatorname{skew}(G_1C_0^{-1})=\operatorname{skew}(B)=A=A_{\mathrm{eff}}\). The naive antisymmetrisation instead keeps the \(C_0\)-weighting: with \(B^\top=S-A\), \[G_1-G_1^\top=BC_0-C_0B^\top=(S+A)C_0-C_0(S-A)=[S,C_0]+\{A,C_0\},\] both terms antisymmetric, so \(\operatorname{skew}(G_1)=\tfrac12[S,C_0]+\tfrac12\{A,C_0\}\). The first term is present already at \(A=0\) and vanishes only when \(S\) and \(C_0\) commute (isotropic \(C_0\)); at the working point it does not (the \(39\%\) quoted). The right-multiplication by \(C_0^{-1}\) removes both this contamination and the \(\{A,C_0\}\) weighting. The identity concerns the one-lag operator of the fitted linear proxy --- the object the estimator forms from data, in which \(A\) is present --- and does not contradict Theorem 4.10, whose strict-Haar population average renders the effective one-lag operator chirality-free; the two describe the fitted proxy and the Haar-averaged population respectively. \(\square\)

\paragraph{4.4.1 Exact metric-ray recovery: a population/oracle result}
\label{exact-metric-ray-recovery}

The background metric in §2.2 is model structure. The next result identifies when an idealized experiment with exact full propagator matrices can recover its inverse ray. This is a population/oracle theorem, not a sample theorem. Write
\[
\mathsf p:=\mathsf g^{-1}:V^*\to V
\]
for the inverse metric.

\bigskip
\textbf{Theorem 4.11.A (observable exact metric-inverse theorem).}
{\itshape Let \(\mathcal T_0,\ldots,\mathcal T_k\in GL(d)\) be exact propagators and define
\[
\mathcal L:\operatorname{Sym}(d)\longrightarrow
\bigoplus_{t=1}^{k}\operatorname{Sym}(d),
\qquad
\mathcal L(Q):=
\bigl(\mathcal T_tQ\mathcal T_t^\top
-\mathcal T_0Q\mathcal T_0^\top\bigr)_{t=1}^{k}.
\tag{4.11a}
\]
Then:

\textup{(i)} \(\ker\mathcal L\) contains a positive-definite \(\mathsf p\) if and only if
\[
\mathcal T_t=\mathcal T_0O_t,
\qquad
O_t\mathsf pO_t^\top=\mathsf p,
\qquad O_t:=\mathcal T_0^{-1}\mathcal T_t.
\tag{4.11b}
\]
In that case the common sandwich
\[
\Gamma:=\mathcal T_t\mathsf p\mathcal T_t^\top
=\mathcal T_0\mathsf p\mathcal T_0^\top
\tag{4.11c}
\]
is independent of \(t\).

\textup{(ii)} Conditional on existence of a positive solution, the inverse-metric ray is unique if and only if \(\dim\ker\mathcal L=1\). The identified object is \(\mathbb R_+\mathsf p\), equivalently \(\mathbb R_+\mathsf g\), not an absolute normalization.

\textup{(iii)} Equivalently, suppose \(\mathcal T_t=BO_t\), with \(O_t\mathsf pO_t^\top=\mathsf p\). The ray is unique if and only if the orthogonalized difference group generated by
\[
\widetilde D_{st}:=
\mathsf p^{-1/2}O_s^{-1}O_t\mathsf p^{1/2}\in O(d)
\tag{4.11d}
\]
acts irreducibly on \(\mathbb R^d\).}

\emph{Proof.}
Part (i) is the common-sandwich equation multiplied by \(\mathcal T_0^{-1}\) and \(\mathcal T_0^{-\top}\). For (ii), if \(Q\in\ker\mathcal L\) is independent of \(\mathsf p\), then \(\mathsf p+\varepsilon Q\succ0\) for all sufficiently small \(|\varepsilon|\), producing more than one positive ray; a one-dimensional kernel containing a positive element contains only that ray. For (iii), whiten a candidate \(Q\) by \(\widetilde Q=\mathsf p^{-1/2}Q\mathsf p^{-1/2}\). The kernel condition is invariance of \(\widetilde Q\) under conjugation by the orthogonal difference group. Irreducibility forces every invariant symmetric matrix to be scalar, because every eigenspace of a non-scalar invariant symmetric matrix would be a proper invariant subspace. Conversely, for a reducible orthogonal group the projector onto a proper invariant subspace is a second symmetric fixed point. \(\blacksquare\)

\bigskip
\textbf{Remark (pointwise dressings preserve feasibility, not the kernel).}
If \(\varepsilon_t\mathsf p\varepsilon_t^\top=\mathsf p\), the dressed propagators \(\mathcal T'_t=\mathcal T_t\varepsilon_t\) retain the feasible true ray:
\[
\mathcal T'_t\mathsf p(\mathcal T'_t)^\top=\Gamma.
\]
In general, \(\mathcal L'\ne\mathcal L\) and their kernels need not agree; uniqueness survives only if the dressed difference group remains irreducible. Indeed, if \(\mathcal T_t=BO_t\), choosing \(\varepsilon_t=O_t^{-1}\) makes every dressed propagator equal to \(B\), so \(\mathcal L'=0\) and \(\ker\mathcal L'=\operatorname{Sym}(d)\). This pointwise right dressing changes the observed propagators. Even coordinate-sign dressings can change nullity: for \(\mathsf p=B=I\), \(P=\left(\begin{smallmatrix}0&0&1\\1&0&0\\0&1&0\end{smallmatrix}\right)\), and \(D=\operatorname{diag}(1,1,-1)\), the family \(\{I,P,P^2,D\}\) has nullity one, whereas dressing only its last element by \(D\) gives \(\{I,P,P^2,I\}\) with nullity two. Pointwise dressing must not be confused with the common factorization gauge \(B\mapsto B\varepsilon\), \(O_t\mapsto\varepsilon^{-1}O_t\), with \(\varepsilon\mathsf p\varepsilon^\top=\mathsf p\), which leaves every \(\mathcal T_t\) unchanged.

\bigskip
\textbf{Scope (oracle versus sample).}
One observed transition gives only \(d\) responses for the \(d^2\) entries of a pointwise propagator. Local-regression estimation of \(\mathcal T_t\) therefore needs local-stationarity, regressor-rank, innovation-orthogonality, and mixing assumptions. Under estimation error, \(\widehat{\mathcal L}\) generally lacks an exact null vector, and its smallest right singular vector need not be positive-definite. A sample theorem would require a propagator-error bound relative to the population singular gap, an SPD-constrained or weighted fit, and calibrated uncertainty. None is established here.

\bigskip
\textbf{Proposition 4.11.B (mixed-moment attenuation and common-Gram repair).}
{\itshape Let
\[
y_{t+1}=BR_ty_t+\eta_t,
\qquad C:=\mathbb E[y_ty_t^\top]\succ0,
\qquad G_1:=\mathbb E[y_{t+1}y_t^\top],
\]
with \(B\in GL(d)\) and \(\mathbb E[\eta_ty_t^\top]=0\). Define \(\mathcal G:=G_1C^{-1}\). Then
\[
\mathcal G=\rho B
\quad\Longleftrightarrow\quad
\mathbb E[R_ty_ty_t^\top]=\rho C.
\tag{4.11e}
\]
The marginal condition \(\mathbb ER_t=\rho I\) is not sufficient by itself; it is sufficient, for example, under \(R_t\perp y_t\), or under \(\mathbb E[R_t\mid y_t]=\rho I\).

Suppose in addition that the exact inverse model supplies the linked-scale pair
\[
\mathsf p\succ0,
\qquad \Gamma=B\mathsf pB^\top,
\]
and that \(\mathcal G=\rho B\). Then
\[
\mathcal G\mathsf p\mathcal G^\top=\rho^2\Gamma,
\qquad
|\rho|^2=\frac1d\operatorname{tr}\!\left(
\Gamma^{-1}\mathcal G\mathsf p\mathcal G^\top\right).
\tag{4.11f}
\]
For \(d=3\), the prefactor is \(1/3\).}

\emph{Proof.}
The cross-moment identity is
\[
G_1=B\,\mathbb E[R_ty_ty_t^\top]
+\mathbb E[\eta_ty_t^\top].
\]
Innovation orthogonality and invertibility of \(B\) make (4.11e) necessary and sufficient. The common-Gram statement follows from
\[
\mathcal G\mathsf p\mathcal G^\top
=\rho^2B\mathsf pB^\top=\rho^2\Gamma;
\]
left-multiplication by \(\Gamma^{-1}\) and normalized trace give (4.11f). \(\blacksquare\)

The scale is linked: \(\mathsf p\mapsto c\mathsf p\) requires \(\Gamma\mapsto c\Gamma\), so \(c\) cancels in (4.11f). Normalizing the rays independently would destroy the identity. Also,
\[
\Gamma^{-1/2}\mathcal G\mathsf p\mathcal G^\top\Gamma^{-1/2}=\rho^2I
\tag{4.11g}
\]
is a necessary model restriction. In a sample, the normalized trace only returns a scalar. One must also assess the proportionality residual, its uncertainty, and a lower confidence bound separating \(|\rho|\) from zero. Since \(A_{\mathsf g}\) is linear,
\[
A_{\mathsf g}(\mathcal G)=\rho A_{\mathsf g}(B).
\tag{4.11h}
\]
Under the exact common-Gram model and \(|\rho|>0\), \(|\rho|\) identifies \(A_{\mathsf g}(B)\) only up to global sign. The scalar-attenuation restriction fixes the common right-\(O(\mathsf p)\) factorization gauge up to the central choices \(\varepsilon=\pm I\); it is not orientation-free. Under the full \(O(\mathsf p)\) model, \((B,\rho)\) and \((-B,-\rho)\) are gauge-equivalent. Choosing the sign requires an external orientation convention.

\bigskip
\textbf{Fresh-Haar boundary.}
For the strict fresh-Haar body recursion, \(R_t\) is conditionally Haar and independent of \(y_t\), so
\[
\mathbb ER_t=0,
\qquad \mathbb E[R_ty_ty_t^\top]=0,
\qquad \mathcal G=0.
\]
Thus \(\rho=0\), and (4.11f) gives no information about \(B\) or \(A_{\mathsf g}(B)\). The nonzero-\(\rho\) repair applies only in a non-Haar or correlated regime supporting both the mixed-moment restriction and common-Gram proportionality.

\bigskip
\textbf{Remark 4.12 (drift bias; Goldie-route bookkeeping).} \emph{The blind-frame drift enters the estimator as \(-\tfrac12\{S,\hat w_{\mathrm{res}}\}\), of the same order as the signal --- the estimator-level counterpart of the irreducible rotation--vibration (Coriolis) coupling that the Eckart frame (Eckart 1935) cannot gauge away. Separately, the Goldie-route split \(D=D_{\mathrm{exp}}+D_\pi\) and the FDT-route lag sum are inequivalent decompositions of the same total; only the totals agree (C21). The explicit-term-only route captures \(55\%\) of the norm at \(31^\circ\); the stationary-law response is the hard part of any closed form.}

We use a quantile/CRN design rather than directly estimating the Goldie constant or \(\partial_\Sigma\log C\) with score-weighted Monte Carlo. The following lemma states the obstruction.

\bigskip
\textbf{Lemma 4.13 (Goldie--Monte-Carlo pathology).} \emph{Let \(R\ge0\) satisfy \(\mathbb P(R>r)\sim Cr^{-\alpha}\), \(C>0\), \(\alpha>2\), independently of \(\mathbf u\sim\mathrm{Unif}(S^2)\) and \(\boldsymbol\eta\sim N(0,\Sigma)\), \(\Sigma\succ0\). For fixed \(S\ne0\), define \(\Phi=\|RS\mathbf u+\boldsymbol\eta\|^\alpha-(R\|S\mathbf u\|)^\alpha\). Then \(\mathbb E|\Phi|<\infty\) but \(\mathbb E[\Phi^2]=\infty\), hence \(\operatorname{Var}(\Phi)=\infty\).} The full proof is in Supplement S.D.4. Exact signed-tail index, stable-domain membership, and trim-bias asymptotics require additional tail-balance arguments and are not claimed here. The reported Goldie level comparison is retained only as a descriptive direction check; its \(\pm0.36\) is not a confidence interval.

\emph{Companion calibration note.} Pre-asymptotic exponent bias is working-point-specific: \(-0.625\pm0.061\) at the positive-definite point and \(-0.718\) at the indefinite point. Legacy calibrations do not transfer. No differential drift was resolved in the reported CRN-paired contrasts (C17).

\begin{center}\rule{0.5\linewidth}{0.5pt}\end{center}

\subsubsection{4.5 Controlled body-frame test: conditional prediction and finite-threshold evidence}\label{body-frame-test-the-effect-is-present-and-predictable}

Sections 4.1--4.2 predict the chiral level shift for a known operator \(B\) and noise frame \(\Sigma\). The controlled simulations test this body-frame prediction directly.

Section 4.2.2 gives the main body-frame result. Under Assumption R, the asymptotic law uses \(D_\infty=\partial_\Sigma\log C\) and \(G_\infty=2H\circ[\Sigma,D_\infty]\). The controlled finite-window comparison instead gives \(\cos\angle(g_{\mathrm{win}},G_{\mathrm{win,pred}})=0.9993\), with \(D_{\mathrm{win}}\) agreeing between the bump and fluctuation--dissipation routes (§4.2.4). The primary statistic is the deep-tail quantile shift \(dq_p\) of (2.6). No differential drift was resolved in the reported CRN-paired contrasts. The hard nulls hold at finite sample; the other simulations are quantitatively consistent at the tested controlled working points.

Blind body-frame recovery does \emph{not} yet follow. Supplying \(Q_t\) creates a controlled frame-oracle experiment, distinct from Theorem 4.11.A's full-propagator oracle. The causal-eigenframe/Kabsch prototype shows qualitative lock loss and Remark 4.12's drift artefact. Its injected-frame design and percentile envelopes establish neither a sampling theorem, confidence interval, nor certified threshold in \(\omega\). It remains only as motivation for another estimator class.

\begin{center}\rule{0.5\linewidth}{0.5pt}\end{center}

\subsubsection{4.6 Empirical sequel}\label{empirical-sequel}

A separate empirical sequel will require a prospectively fixed design, complete provenance, and an independent sample.

\begin{center}\rule{0.5\linewidth}{0.5pt}\end{center}

\subsubsection{4.7 What §4 establishes; open problems → companion work}\label{what-4-establishes-open-problems-companion-work}

\bigskip
\textbf{Established} (collected in §5.1). The background metric is necessary model structure. The low-order odd scalar channel pairs with metric-relative covariance misalignment (Propositions 2.1--2.2, Theorem 4.0). The polar-twist mechanism is exact (Theorem 4.1). Under Assumption R, the first-order response and confinement law hold; finite-threshold simulations are consistent with them at the working points. In the strict-Haar model, every lagged second-order cross-moment is transpose-even, while a fourth radial moment retains an odd term (Theorem 4.10, Corollary 4.10.D, Proposition 4.10.B).

\bigskip
\textbf{Resolved negative (not pursued in the companion work).} No singular resonance or new critical mechanism is established. Within the nonzero-eigenvalue chamber, the apparent pair-sum pole multiplies a closed twist channel, and the near-pairing sweep remains regular. The indefinite chamber still differs structurally from the SPD interior: a regular direct channel carries its mixed-sign response, whose closed gradient theory remains Open Problem 2.

\bigskip
\textbf{Open problems.} 1. \textbf{Conditional-score profile and \(\varphi\).} Close the stationary-law response \(D_\pi\) (Goldie route) or the \(k\ge2\) tilted sum (FDT route), instead of the quantitatively excluded explicit-term-only approximation (Remark 4.12). 2. \textbf{Mixed-sector gradient theory.} Develop an independent theory for the symmetry-free mixed-sign sector of \(\operatorname{Skew}\nabla_B\log C\) (Theorem 4.3). 3. \textbf{Blind body-frame recovery.} Replace the causal-eigenframe prototype with a state-space or rotation-increment estimator and a genuine sampling theorem; no recovery frontier is certified. 4. \textbf{Sample metric-ray inference.} From Theorem 4.11.A, derive propagator-error and singular-gap bounds, an SPD-constrained fit, a proportionality test for (4.11g), and uncertainty for \(|\rho|\).

\begin{center}\rule{0.5\linewidth}{0.5pt}\end{center}

\subsection{5. Discussion}\label{discussion}

\subsubsection{5.1 Summary}\label{summary}

This paper characterises transpose-odd operator content relative to an explicit positive metric. Positive mirrors form a five-dimensional metric-ray moduli space, so no \(GL(3)\)-natural metric-free positive transpose exists (Proposition 2.1). At fixed metric, \(B\)-only odd polynomial scalars first occur at degree six. Independent covariance opens a complete low-order pairing family with common zero criterion \([S_{\mathsf g},\Sigma\mathsf g]=0\) (Proposition 2.2, Theorem 4.0). In the fresh-Haar Kesten model, the tail exponent is transpose-blind, but the level need not be. The exact polar twist gives an unconditional response mechanism. A nonzero response is established conditionally, under the displayed regularity hypotheses, through the first-order law and fluctuation--dissipation representation, together with finite-threshold numerics. Strict-Haar lagged second-order moments do not identify the odd operator, while a fourth radial moment retains an odd term. No finite scalar-weighted radial-moment combination uniformly cancels the even quadratic form while retaining the odd covector. The exact metric-ray inverse and common-Gram repair remain population/oracle results only.

\subsubsection{5.2 Scope and limitations}\label{scope-and-limitations}

Three limits bound the results. First, the background metric is explicit model structure, fixed as \(\mathsf g=I\) in the canonical chart and independent of \(\Sigma\). Second, the exact confinement kernel and working-point simulations are three-dimensional, although the polar-twist mechanism is dimension-agnostic. Third, the first-order level theory requires Assumption R, and the fluctuation--dissipation representation requires Assumption UI″. The proved fixed-lag nullity, lag-window escape, and Appendix C density bound do not prove these full regularity assumptions. Noisy use of Theorem 4.11.A would require propagator-error control, a population singular gap, SPD inference, and uncertainty calibration.

\subsubsection{5.3 Outlook}\label{outlook}

Three tracks remain. First, Proposition 4.10.B gives a linear-fit-free fourth-moment odd channel, while Proposition 4.10.C explains why operator-uniform radial combinations cannot isolate its sign from even magnitude. Second, blind co-moving-frame recovery needs a state-space estimator and sampling theorem; the prototype gives no certified threshold. Third, Theorem 4.11.A motivates SPD-constrained metric-ray inference for noisy local propagators, with singular-gap and common-Gram proportionality tests. Other problems are analytic closure of \(\varphi\) and correlated multi-step models in which operator asymmetry might reach the exponent.

\subsection{Appendix A. Checkpoint and status crosswalk}\label{appendix-a.-checkpoint-table-c1c31}

Checkpoint identifiers are diagnostics, not proofs. Supplement S.0 is the controlling result-status ledger; Supplements S.E.2--S.E.3 contain the checkpoint index and reproduction audit.

\subsection{Appendix B. Assumption and proof crosswalk}\label{appendix-b.-assumptions-proofs-and-the-conversion-identity}

Assumptions R and UI\(^{\prime\prime}\) precede their first uses in §§4.1.2 and 4.2.4. Supplement S.B.8 gives their scope and justification audit, expanded proofs, fixed-\(p\)/fixed-\(u\) conversion, and \(d=4\) confinement formula. Supplement S.C.1 gives the fluctuation--dissipation and coupling audit.

\subsection{Appendix C. A local tail-density bound (Theorem LR)}\label{appendix-c.-a-local-tail-density-bound-theorem-lr}

Lemma 4.8.1(b)'s rate uses the following local-density bound for the stationary radial law. It shows that one-step singularities in the multiplicative gain density self-heal in the stationary tail. Standard Kesten--Goldie theory controls the tail probability \(\mathbb P(R>z)\). Supplement S.C.2 supplies the density-level statement and complete proof, with no companion-work dependency.

\bigskip
\textbf{Setting.} \(r_{t+1}=\|r_tB\hat{\mathbf u}_t+\boldsymbol\eta_t\|\) in \(\mathbb R^3\); \(\hat{\mathbf u}_t\sim\mathrm{Unif}(S^2)\), \(\boldsymbol\eta_t\sim\mathcal N(0,\Sigma)\), \(\Sigma\succ0\), i.i.d. and independent; \(B\) invertible with singular values \(\sigma_1\ge\sigma_2\ge\sigma_3>0\), not all equal; \(R\) the stationary radial law. The standing §3 inputs are the root equation and negative drift for \(M:=\|B\hat{\mathbf u}\|\), namely \(\mathbb E M^{\alpha_\star}=1\), \(\gamma:=\mathbb E\log M<0\) (H), and the tail asymptotic \(\bar F(z):=\mathbb P(R>z)=Cz^{-\alpha_\star}(1+\varepsilon(z))\), \(C>0\), \(\varepsilon(z)\to0\) (T).

\bigskip
\textbf{Theorem LR.} \emph{There are \(\delta_0>0\), \(z_0\), \(C_\sharp\), depending only on \((B,\Sigma,\alpha_\star,C,\varepsilon(\cdot))\), such that for \(z\ge z_0\),} \[
\sup_{y\in[z/2,2z]}f_R(y)\le C_\sharp\,\frac{\bar F(z)}z,\qquad \mathbb P(|R-y|\le u)\le C_\sharp\,u\,\frac{\bar F(y)}y\quad(0<u\le\delta_0y,\ y\ge z_0).
\]

\emph{Role and proof.} Lemma 4.8.1(b) uses only the second inequality, with \(y=z\) and \(u=D_k(z)\); the first is the accompanying density sup-bound. Supplement S.C.2 gives the six-step proof: Gaussian radial kernel, angular-gain law, stationary-shell pullback, multiscale contraction, fine \(|V|\)-shells, and layer-cake assembly. Its inputs are the displayed \(d=3\) setting, full-rank Gaussian noise, invertibility of \(B\), and the Section 3 tail asymptotic (T). The result is a sup-bound, not a pointwise density asymptotic. It proves neither Assumption UI\(^{\prime\prime}\)(ii) nor convergence \(\alpha_{\mathrm{loc}}\to\alpha_\star\).

\bigskip
\textbf{Additional citation context for §4} (beyond the Kesten/Goldie/Buraczewski--Damek--Mikosch/Bougerol--Picard set cited in §2--§3): Procesi (1976, Theorem 7.1) for orthogonal trace-word invariants; Glynn (1990) and L'Ecuyer (1990) for the likelihood-ratio/score gradient (Theorem 4.8); Kubo (1966) and Marconi--Puglisi--Rondoni--Vulpiani (2008) for the fluctuation--dissipation framing; Kenney--Laub (1991) and Higham (2008, Ch.~8) for the polar-factor Fréchet derivative; Glynn--Olvera-Cravioto (2019), Meyn--Tweedie, and Kloeckner (2022) for stationary-measure sensitivity; Straumann--Mikosch (2006) for analogous SRE differentiability; Hillier (2001) for quadratic-form-on-sphere densities; and Collamore--Vidyashankar (2013) for tilted/regeneration machinery. Troude--Sornette (2025, 2026) is cited only as heuristic and numerical context; no approximation from those papers is used.

\emph{The controlling result-status ledger and reproducibility crosswalk are Supplement S.0 and S.E.2--S.E.3.}

\subsection*{Data and code availability}
No new empirical dataset is analysed. The Supplement provides the mathematical derivations and the descriptions of the numerical checkpoints. The code used for the numerical checkpoints is available from the author upon request.

\subsection*{Declaration of competing interests}
The author declares no known competing financial interests or personal relationships that could have appeared to influence the work reported in this paper.

\subsection*{Declaration of generative AI and AI-assisted technologies}
All theorems and proofs* are the author's, and the author takes full responsibility for them. The text was translated from Turkish into English with Claude Opus 4.8 (Anthropic), which was also used for language editing. AI tools were used to write and run the numerical checks; the author reviewed their outputs. Claude Fable 5 (Anthropic) and GPT-5.6 Sol (OpenAI) were used for cross-checking. *Following a warning raised by Sol during cross-checking, the author, working together with Sol, added the independent-copy coupling proof of Lemma 4.8.1 and Corollary 4.8.2, and, on Sol's recommendation, decided to add Appendix C. Both contributions were directed by the author, and responsibility for them rests with the author.

\bigskip
\textbf{References.}

\begin{itemize}
\tightlist
\item
  Alsmeyer, G., Brofferio, S., Buraczewski, D. (2023). Asymptotically linear iterated function systems on the real line. \emph{Annals of Applied Probability} \textbf{33}(1), 161--199.
\item
  Ang, A., Chen, J. (2002). Asymmetric correlations of equity portfolios. \emph{Journal of Financial Economics} \textbf{63}, 443--494.
\item
  Bougerol, P., Picard, N. (1992). Strict stationarity of generalized autoregressive processes. \emph{Annals of Probability} \textbf{20}(4), 1714--1730.
\item
  Buraczewski, D., Damek, E. (2017). A simple proof of heavy-tail estimates for affine-type Lipschitz recursions. \emph{Stochastic Processes and their Applications} \textbf{127}, 657--668. \url{https://arxiv.org/abs/1604.06943}.
\item
  Buraczewski, D., Damek, E., Mikosch, T. (2016). \emph{Stochastic Models with Power-Law Tails: The Equation X = AX + B}. Springer.
\item
  Collamore, J. F., Vidyashankar, A. N. (2013). Tail estimates for stochastic fixed-point equations via nonlinear renewal theory. \emph{Stochastic Processes and their Applications} \textbf{123}, 3378--3429.
\item
  Diaconis, P., Freedman, D. (1999). Iterated random functions. \emph{SIAM Review} \textbf{41}(1), 45--76.
\item
  Eckart, C. (1935). Some studies concerning rotating axes and polyatomic molecules. \emph{Physical Review} \textbf{47}(7), 552--558.
\item
  Elton, J. H. (1990). A multiplicative ergodic theorem for Lipschitz maps. \emph{Stochastic Processes and their Applications} \textbf{34}(1), 39--47.
\item
  Glynn, P. W. (1990). Likelihood ratio gradient estimation for stochastic systems. \emph{Communications of the ACM} \textbf{33}(10), 75--84.
\item
  Glynn, P. W., Olvera-Cravioto, M. (2019). Likelihood ratio gradient estimation for steady-state parameters. \emph{Stochastic Systems} \textbf{9}(2), 83--100.
\item
  Goldie, C. M. (1991). Implicit renewal theory and tails of solutions of random equations. \emph{Annals of Applied Probability} \textbf{1}(1), 126--166.
\item
  Higham, N. J. (2008). \emph{Functions of Matrices: Theory and Computation}, Ch. 8. SIAM.
\item
  Hillier, G. (2001). The density of a quadratic form in a vector uniformly distributed on the \(n\)-sphere. \emph{Econometric Theory} \textbf{17}(1), 1--28.
\item
  Kenney, C., Laub, A. J. (1991). Polar decomposition and matrix sign-function condition estimates. \emph{SIAM Journal on Scientific and Statistical Computing} \textbf{12}, 488--504.
\item
  Kesten, H. (1973). Random difference equations and renewal theory for products of random matrices. \emph{Acta Mathematica} \textbf{131}, 207--248.
\item
  Kloeckner, B. R. (2022). Optimal transportation and stationary measures for iterated function systems. \emph{Mathematical Proceedings of the Cambridge Philosophical Society} \textbf{173}, 163--187.
\item
  Kubo, R. (1966). The fluctuation--dissipation theorem. \emph{Reports on Progress in Physics} \textbf{29}(1), 255--284.
\item
  L'Ecuyer, P. (1990). A unified view of the IPA, SF, and LR gradient-estimation techniques. \emph{Management Science} \textbf{36}(11), 1364--1383.
\item
  Longin, F., Solnik, B. (2001). Extreme correlation of international equity markets. \emph{Journal of Finance} \textbf{56}(2), 649--676.
\item
  Marconi, U. M. B., Puglisi, A., Rondoni, L., Vulpiani, A. (2008). Fluctuation--dissipation: response theory in statistical physics. \emph{Physics Reports} \textbf{461}(4--6), 111--195.
\item
  Meyn, S. P., Tweedie, R. L. (2009). \emph{Markov Chains and Stochastic Stability}, 2nd ed.~Cambridge University Press.
\item
  Patton, A. J. (2006). Modelling asymmetric exchange-rate dependence. \emph{International Economic Review} \textbf{47}(2), 527--556.
\item
  Procesi, C. (1976). The invariant theory of \(n\times n\) matrices. \emph{Advances in Mathematics} \textbf{19}(3), 306--381. The orthogonal trace-word result used here is Theorem 7.1. \url{https://doi.org/10.1016/0001-8708(76)90027-X}.
\item
  Smith, G. F. (1971). On isotropic functions of symmetric tensors, skew-symmetric tensors and vectors. \emph{International Journal of Engineering Science} \textbf{9}(10), 899--916. \url{https://doi.org/10.1016/0020-7225(71)90023-1}.
\item
  Spencer, A. J. M., Rivlin, R. S. (1962). Isotropic integrity bases for vectors and second-order tensors. \emph{Archive for Rational Mechanics and Analysis} \textbf{9}, 45--63. \url{https://doi.org/10.1007/BF00253332}.
\item
  Straumann, D., Mikosch, T. (2006). Quasi-maximum-likelihood estimation in conditionally heteroscedastic time series: a stochastic recurrence equations approach. \emph{Annals of Statistics} \textbf{34}(5), 2449--2495.
\item
  Troude, V., Sornette, D. (2025). Eigenvector geometry as a new route to criticality in random multiplicative systems. arXiv:2510.21755.
\item
  Troude, V., Sornette, D. (2026). Non-normal eigenvector amplification in multidimensional Kesten processes. \emph{Physical Review Research} \textbf{8}(3), 033050. \url{https://doi.org/10.1103/ztdv-7t5j}. arXiv:2510.11763.
\item
  Vigdorovich, I., Foysi, H. (2015). Simultaneous invariants of strain and rotation rate tensors and their admitted region. \emph{Advances in Mathematical Physics} \textbf{2015}, Article ID 147125. \url{https://doi.org/10.1155/2015/147125}.
\item
  Zheng, Q.-S. (1993). On the representations for isotropic vector-valued, symmetric tensor-valued and skew-symmetric tensor-valued functions. \emph{International Journal of Engineering Science} \textbf{31}(7), 1013--1024. \url{https://doi.org/10.1016/0020-7225(93)90109-8}.
\item
  Zheng, Q.-S. (1994). Theory of representations for tensor functions---a unified invariant approach to constitutive equations. \emph{Applied Mechanics Reviews} \textbf{47}(11), 545--587. \url{https://doi.org/10.1115/1.3111066}.
\end{itemize}
\end{document}